\documentclass[11pt]{article}
\usepackage{epsfig}
\usepackage{amsmath}
\usepackage{amssymb}
\usepackage{graphicx}
\usepackage[latin1]{inputenc}
\def\BibTeX{{\rm B\kern-.05em{\sc i\kern-.025em b}\kern-.08em
    T\kern-.1667em\lower.7ex\hbox{E}\kern-.125emX}}

\def\limiten{\renewcommand{\arraystretch}{0.5}
\begin{array}[t]{c}\stackrel{}{\longrightarrow} \\
{\scriptstyle n\rightarrow
\infty}\end{array}\renewcommand{\arraystretch}{1}}

\def\limitep{\renewcommand{\arraystretch}{0.5}
\begin{array}[t]{c}\stackrel{}{\longrightarrow} \\
{\scriptstyle \lambda \rightarrow
\infty}\end{array}\renewcommand{\arraystretch}{1}}

\def\limiteloin{\renewcommand{\arraystretch}{0.5}
\begin{array}[t]{c}\stackrel{{\cal D}}{\longrightarrow} \\
{\scriptstyle n\rightarrow
\infty}\end{array}\renewcommand{\arraystretch}{1}}

\newcommand{\be}{\begin{equation}}
\newcommand{\ee}{\end{equation}}
\newcommand{\bd}{\begin{displaymath}}
\newcommand{\ed}{\end{displaymath}}

\newcommand{\ba}{\begin{eqnarray}}
\newcommand{\ea}{\end{eqnarray}}
\newcommand{\ban}{\begin{eqnarray*}}
\newcommand{\ean}{\end{eqnarray*}}

\newcommand{\LL} {I\!\!L}
\newcommand{\C} {I\!\!\!\!C}
\newcommand{\R} {I\!\!R}
\newcommand{\E} {I\!\! E}
\newcommand{\N} {I\!\! N}
\newcommand{\Z} {{\bf Z}}
\newcommand{\cov} {\mbox{cov}\hspace{0.9mm}}
\newcommand{\var} {\mbox{var}\hspace{0.9mm}}

\renewcommand{\arraystretch}{.8}

\renewcommand{\Box}{\hfill\rule{0.25cm}{0.25cm}} 

\newtheorem{Prop}{Proposition}
\newtheorem{lem}{Lemma}
\newtheorem{Theo}{Theorem}
\newtheorem{cor}{Corollary}
\newtheorem{rem}{Remark}

\newenvironment{dem}{\  {\bf Proof. }}
{\Box\par\medskip\noindent}
\newcommand{\ds}{\displaystyle}
\def\1{{\bf 1}}
\def\Pr{I\mskip-7muP}

\begin{document}


\title{A nonparametric estimator of the spectral density of a continuous-time Gaussian process observed at random times}
\author{JEAN-MARC~BARDET${}^*$ and PIERRE, R.~BERTRAND${}^{**, \dag}$\\
~\\
${}^*$  {\it SAMOS, Université Paris 1,} \\
${}^{**}$ {\it INRIA Saclay and Laboratoire de
Mathématiques, Université Clermont-Ferrand.} }
\pagestyle{myheadings}

\maketitle

\begin{abstract}
In numerous applications data are observed at random times and an estimated graph of the spectral density may be relevant for characterizing and explaining phenomena. By using a wavelet analysis, one derives a nonparametric estimator
of the spectral density of a Gaussian process with stationary
increments (or a stationary Gaussian process) from the observation of one path  at random discrete times. For every positive frequency, this estimator is proved to satisfy a central limit theorem with a convergence rate depending on the roughness of the process and the moment of random durations between successive observations.
In the case of stationary Gaussian processes, one can compare this estimator with estimators based on the empirical periodogram. Both estimators reach the same optimal rate of convergence, but the estimator based on wavelet analysis converges for a different class of random times.
Simulation examples and application to biological data are also provided.
\end{abstract}
\section{Introduction}
In biology, finance, internet traffic, oceanography, civil engineering, etc..,
detrended data are often modeled by centered Gaussian processes observed at random times.
Under some stationarity assumption (or by assuming the stationarity of increments), such processes are characterized by their spectral density function.
Roughly speaking the spectral density function corresponds to the Fourier transform of the covariance of the process (or its increments), and in  the sequel, it will be denoted by $f(\xi)$, where $\xi$ is the frequency.

Before going further, let us give a detailed example: the heart rate variability. Cardiologists are interested in the behavior of its  spectral density,
usually on both frequency bands 
$(\omega_1,\,\omega_{2})=(0.04\,Hz,\,0.15\, Hz)$ and $(\omega_2,\,\omega_{3})=(0.15\, Hz,\,0.5\, Hz)$ 
corresponding  to the orthosympathetic nervous system and the parasympathetic one, respectively, see Task force of the European Soc. Cardiology and the North American Society
of  Pacing and Electrophysiology (1996). The spectral density follows different power laws on the different frequency bands, {\em i.e.} $f(\xi) = \sigma_i\,|\xi|^{-\beta_i}$ when $\xi\in(\omega_i,\omega_{i+1})$.
Finally, according to the type of activity or the period of the day, we notice variations of these parameters, see Section \ref{cardio}.

As shown by the previous example, the spectral density contains relevant information.  The phenomenon of interest is known through the observation of a Gaussian stationary process at random times.
From a statistical point of view, there are two different situations regarding estimation of the spectral density. On the one hand,  for Gaussian processes with stationary increments (which include fractal processes), most of the statistical studies  concern  the behavior of the spectral density when $|\xi| \to\infty$ or in the neighborhood of $0$, with a regularly spaced sampling scheme, see Dahlhaus (1989), Gloter and Hoffmann (2007), Moulines {\em et al.} (2007), Begyn (2005) or the book edited by Doukhan {\it et al.} (2003).
On the other hand,  the estimation of the spectral density of stationary Gaussian processes is a classical problem, which has been studied by using periodogram methods, see Shapiro and Silverman (1960), Parzen (1983), Masry (1978a-b) or Lii and Masry (1994). A precise condition to avoid {\em aliasing} has been given in Masry (1978a-b), see  Section \ref{setting:problem}. Unfortunately, Masry's Conditions are not satisfied in the heartbeat time series example. In this case,  the durations between two observation times are bounded since they correspond to the time  between two successive heartbeats. Thus, we have introduced a different modelling of observation times, see Section \ref{setting:problem}.

Moreover, we prefer to use wavelet analysis rather than empirical periodograms. This approach was introduced for processes with stationary increments (in the particular case of fractional Brownian motion) by Flandrin (1992), see also Abry {\it et al.} (2003).
From wavelet analysis, we derive a nonparametric estimator of the spectral density.
Then this estimator is proved to satisfy a central limit theorem (CLT in the sequel).
Numerical applications show the good accuracy of this estimator in the case of processes having stationary increments, e.g.
fractional Brownian motion, as well as  in the case of a stationary processes, e.g. the Ornstein-Uhlenbeck process.

The rest of the paper is organized as follows:
Section~\ref{setting:problem} gives a description of the problem.
Section~\ref{maths} is devoted to  wavelet analysis 
and the CLT satisfied by the estimator of the spectral density $f$. This estimator is applied to generated and real data in Section \ref{numeric}. Appendix contains the proofs.
\section{Description of the problem}\label{setting:problem}
Consider first a Gaussian process $X=\{X(t), t\in\R\}$ with zero
mean and stationary increments. Results can be however extended to the case where a polynomial trend is added to such processes. Therefore $X$ can be written following a
harmonizable representation, see Yaglom (1958) or Cram\`er \& Leadbetter (1967). We adopt a more recent notation as in
Bonami \& Estrade (2003), thus
\be \label{repr:harmonizable}
X(t)= \int_{\R}\big(
e^{it\xi} - 1\big) f^{1/2}(\xi) \, dW(\xi),~~~~\mbox{for
all}~~t \in \R, \ee where
\begin{itemize}
\item the complex isotropic random measure
$dW$ satisfies $dW = dW1 + i\, dW2$ with
$dW1$ and $dW2$ two independent real-valued Brownian measures
(see more details on this part in section 7.2.2 of Samorodnitsky \& Taqqu (1994). Therefore when $g=g_1 + i\,  g_2$ where $g_1$ and $g_2$ are  respectively even and odd real-valued functions such that $\int_{\R} (g_i^2(x))dx<\infty$ ($i=1,2$) then $\E \big [ \big (\int_{\R} g(\xi)\, dW(\xi) \big)^2 \big]=\int_{\R} |g(x)|^2dx<\infty$. Moreover, if $h=h_1 + i \, h_2$ where $h_1$ and $h_2$ are also respectively even and odd real-valued functions such that $\int_{\R} (h_i^2(x))dx<\infty$ ($i=1,2$), then
\be
\label{isometry:IS}
\E \big[\big (\int_{\R} g(\xi)\, dW(\xi) \big)\big (\int_{\R} h(\xi)\, dW(\xi) \big)\big] = \int_{\R} g(x) \overline{h(x)}\, dx.
\ee

\item the function $f$ is called the spectral density of $X$ and is a non-negative even function such that
\be \label{cond:f} \int_{\R} \big(1\wedge
|\xi|^2\big)\, f (\xi) \, d \xi < \infty. \ee
In the sequel, $f$ will be supposed to
satisfy also Assumption F($H$) defined in subsection \ref{nonpar}, but conditions are weak and the class of processes that can be considered is general. 
\end{itemize}
As a particular case, if $X$ is a
stationary processes, one still denote by $f$ the spectral density such that
\begin{eqnarray} \label{stat} X(t)= \int_{\R}
e^{it\xi}\, f^{1/2}(\xi) \, dW(\xi)~~~\mbox{for
all}~~t \in \R,
\end{eqnarray}
where $f$ is still a Borelian positive even function, but satisfying the stronger condition
\be \label{cond:f:stat} \int_{\R}  f (\xi) \, d \xi < \infty. \ee
Even if the definitions are different, $f$ denotes as well the spectral density of a process having
stationary increments or a stationary process. Indeed, in the sequel we consider wavelet coefficients of $X$  which have the same expression with respect to $f$ for both models (\ref{repr:harmonizable}) and (\ref{stat}), see more details 
in Proposition \ref{pte:d(a)}.
Define also the $\sigma$-algebra ${\cal F}_X$
generated by the process $X$, {\it i.e.} \ba \label{tribu} {\cal
F}_X := \sigma \big \{X(t),\, t \in \R \big \}.
\ea
A path of such a process
$X$ on the interval $[0, T_{n}]$ at discrete times $t^{(n)}_i$ for
$i=0,1,\dots, n$ is observed, {\it i.e.}
\begin{eqnarray*}
\label{observe} \big
(X(t^{(n)}_0),X(t^{(n)}_1),\ldots,X(t^{(n)}_n)\big )~~\mbox{is
known, with}~~0=t^{(n)}_0<t^{(n)}_1<\dots<t^{(n)}_{n}=T_n.
\end{eqnarray*}
A unified frame of irregular observed times, grouping deterministic and stochastic
ones, will be considered. First let us assume that there exist a sequence of positive real numbers $(\delta_n)_{n \in \N}$ and
a sequence of random variables  (r.v. in the sequel)
$(L_k)_{k\in \N }$ (which could be deterministic real numbers) such that
\ba \label{lag} \forall k \in
\{0,1,\ldots,n-1\},~~~  t^{(n)}_{k+1}-t^{(n)}_{k} &:=&
\delta_n \, L_k,\quad \mbox{and}~~\delta_n \limiten 0.\vspace{-0.5cm}
\ea
\noindent
For $Z$ a r.v. and $\alpha\in (0,\infty)$, denote $\| Z\|_\alpha:=\big ( \E (|Z|^\alpha) \big ) ^{1/\alpha}$
if $\E (|Z|^\alpha) <\infty$.
Now, assume that there exists $s \in [1,\infty)$ such that
\vspace{0.2cm}
\newline {\bf Assumption S($s$)} {\em $(L_k)_{k \in \N}$ is a sequel
of 
positive r.v. such that:
there exist $0<m_1\leq M_1$ and $M_s<\infty$ satisfying
$$
m_1\le \E L_k=\|L_{k}\|_1\le M_1\quad \mbox{and} \quad \|L_{k}\|_s \leq M_s,\quad \mbox{for all $k\in \N$}.
$$
}
\noindent Then we can also define:\\
\newline {\bf Assumption S($\infty$)} {\em $(L_k)_{k \in \N}$ is a sequence
of 
positive r.v. satisfying Assumption S($s$) for every $s\in \N$.}\\
~\\
For instance, it is clear that if $(L_k)_{k \in \N}$ is a sequence of exponential or bounded r.v., then Assumption S($\infty$) is satisfied.
Now, $\ds T_n=\delta_n \big ( L_0+\ldots+L_{n-1} \big )$ and under Assumption S($s$) for any $s\geq 1$,
\be
\label{Tn:eq:ndelta}
m_1\times (n\delta_n)\le \E\big(T_n\big)\le M_1\times (n\delta_n).
\ee
This point will be extensively used in the sequel to replace the asymptotic $\E\big(T_n\big) \to \infty $ by $n\delta_n\to \infty$.
\subsection*{Comments on the modelling of observation times}
Assumption S($s$) on the observation times may seem slightly unusual. This leads to the
following comments:
\begin{enumerate}
\item  
Generally, for processes observed at random times (see for instance Lii and Masry, 1994), the duration between observation times $\tau_k=(t_{k+1}-t_{k})$ are  random variables not depending on the data length. In a sense, the asymptotic behavior only concerns the length of observation time $T_n$ (which is necessary to estimate the spectral density at low frequencies). But the lag between two successive random observation times have to be, sufficiently often, very small to allow an estimation of the spectral density for high frequencies. Hence, observation times have to satisfy a strong condition and are typically of a Poissonian type.
\item
In our modelling there are two asymptotic behaviors: the length of the observation time $T_n$ converges to infinity and the mesh $\delta_n$ converges to $0$. The first one is standard up to the slight difference that $T_n$ can be random. This assumption is justified by numerous applications; for instance, the duration of a marathon is clearly random. Thus we have to replace the first asymptotic by $\E\big(T_n\big) \to \infty$.
The second one is less standard but corresponds to applications. We have followed and transposed the idea of round-off introduced, to our best knowledge, by Delattre and Jacod (1997) and currently used today, see for example Robert and Rosenbaum (2008).
In this setting, the time is continuous but round-off with a precision $\delta_n$. Then, the duration between observation times $\big(t_{k+1}^{(n)}-t_{k}^{(n)}\big)$ are the mesh
$\delta_n$ multiplied by positive random variables $L_k$.
Eventually, we do not assume that the r.v. $L_k$ are independent nor identically distributed.
 \item Our choice which is also relative to numerous application cases (see the example of heart rate variability below) has been to provide a spectral density estimation under very weak conditions on the observation times. Typically our results remain valid even for regular sampling, which is not the case under Masry's conditions defined below.
\item In applications, signals are observed at discrete times which are mostly irregularly spaced and random. This type of observations can be met in medicine, physics, mechanics, oceanography,... In these cases observation times depend on the measuring instrument, therefore of a random independent from that of the process $X$. In this context, the hypothesis of independence of durations $(L_k)_{k \in \N}$ and ${\cal F}_X$  is completely realistic. The only case where this assumption seems restrictive concerns financial data. However it is until this day always made, see for instance, Hayashi and Yoshida (2005) or A\"{\i}t-Sahalia and Mykland (2008).
\end{enumerate}
\subsection*{Estimation of the spectral density, state of the art}
To our knowledge, the estimation of the spectral density of a Gaussian process with stationary increments on finite bands of frequencies, from observation at discrete times, is a new problem. Recall that the spectral density $f(\xi)= C\, |\xi|^{-(2H+1)}$ corresponds to a fractional Brownian motion (fBm in the sequel) with Hurst index $H$. However, most of the statistical studies devoted to the fBm or its generalizations concern estimation of the local regularity parameter (linked to the behavior of the spectral density at $\infty$) or the long memory parameter (linked to the behavior of the spectral density in the neighborhood of $0$). The estimation of the spectral density of stationary Gaussian processes is a classical problem corresponding to numerous practical applications, see Shapiro and Silverman (1960) or Parzen (1983). The used methods are based on the periodogram defined by $ I_T(\xi) =  (2 \pi T)^{-1} \Big |\int_0^T e^{-i\xi t} X(t)\, dt \Big |^2$. However, if $(X_{t_1},\cdots,X_{t_N})$ is known, with  regularly spaced observation times $t_i= i\Delta$ and  $T=t_N=N\Delta $, then $  \lim_{T\to \infty} \E I_T(\xi) = f(\xi)$ but $ \lim_{N\to \infty} \E J_N(\xi) = \sum_{k\in\Z}f(\xi+2 k \pi \Delta^{-1})$ where $J_N$ is the empirical periodogram, that is $ J_N(\xi) := (2 \pi N\Delta)^{-1} \Big|\Delta\sum_{k=1}^N e^{-i\xi k \Delta} X(k \Delta) \Big|^2$. Such a phenomenon is called   {\em aliasing}.
To avoid {\em aliasing}, random sampling is chosen and then the empirical periodogram becomes asymptotically unbiased. By using a spectral window an estimator of the spectral density can be deduced and it satisfies a central limit theorem (CLT) with a rate of convergence $T^{-2/5}$, see Masry (1978a-b) or Lii and Masry (1994). These results are obtained for random sampling satisfying very specific conditions that we will call in the sequel:
\begin{description}
\item[Masry's conditions:] the process of  observation times  $\big(t_k\big)_k$ is a stationary, orderly point process independent of $X$, with known mean rate $\beta$ and covariance density $c( u ) $ and satisfies the condition $\beta^2+ c(u)>0$ a.e., where $N(\cdot)$ is the associate counting process, $\beta= \E\big[ N\big((0,1)\big)\big]$ and $c(u)$ its covariance density function (Masry, 1978, Cor. 1.1, p. 320).
\end{description}
When the trajectory is not sampled but observed at random times not chosen by the experimenter, a first step before the estimation of the spectral density is to check that the family $(t_i)$ satisfies Masry's conditions and for this it is necessary to estimate the mean rate $\beta$ and the covariance density function $c(u)$.
\subsection*{Wavelet based estimators}
Wavelet analysis was already used to estimate the parametric behavior of a power law spectral density when $\log|\xi| \to \infty$ or $\log|\xi| \to -\infty$ in a time series (with regularly spaced observation times). In the sequel, we will show that the wavelet analysis is also an interesting tool to estimate the spectral density for Gaussian processes having stationary increments (or stationary Gaussian processes) when a path is observed at random times. Let us underline that the wavelet analysis in Abry {\em et al.} (2003) is based on the sample variance of wavelet coefficients and thus is different from that proposed by Lehr and Lii (1997) or  Goa {\em et al.} (2002) who respectively consider the wavelet decomposition of the estimator derived from the empirical periodogram and the periodogram of the Haar wavelet transform of the process. In both these last cases, discrete time observations are supposed to satisfy Masry's conditions to avoid aliasing. \\
~\\
We consider a non-parametric estimator of the spectral density based on a sample variance of wavelet coefficients. There are two main differences with the approach of Flandrin (1992) or Abry {\em et al.} (2003). Firstly, the definition of ``empirical'' wavelet coefficients, see  (\ref{def:eab}),
is adapted for non-regular observation times. Then a general CLT for sample variance of such ``empirical'' wavelet coefficients is established (see theorem \ref{Theo:TCL:discret}) and a CLT for a semiparametric estimator of the spectral density can be deduced for a large class of fractional processes.
Secondly, one considers a sequel of mother wavelets $\psi_{\lambda}$
in a way that enables the convergence, as $\lambda\to \infty$ of $|\widehat{\psi}_\lambda|^2$ to a Dirac mass concentrated at the frequency $\xi=1$.
Then a CLT for a nonparametric estimator of the spectral density is derived (see Proposition \ref{cvg:Ip}). For observation times satisfying Assumption $S(s)$ with $s> 2$, the supremum of the convergence rate of this last CLT is $T_n^{-2/5}$. This is the same convergence rate as for the periodogram based estimator one (see for instance Lii and Masry, 1994), but for a class of observation times clearly more general than the Masry's one (see for example  Lii and Masry, 1994). Indeed,  our assumptions on observation times allow non-stationary or regularly spaced times, for Gaussian stationary processes
and also for Gaussian processes having stationary increments (like fBm). However, a relation between $T_n$ and $\delta_n$ is required (see condition \ref{Cond:delta:T} below). This condition depends on the regularity of the trajectory and the variability $s$ of observation times. Therefore, in terms of the number $n$ of observations, the convergence rate of our estimator $\ds\widehat f^{(\lambda_n)}_n(\xi)$ is slower than $n^{-2/5}$. \\
~\\
Finally, let us add two comments on the choice and the advantage of wavelet based estimators. Firstly, our method plainly uses the time-frequency localization of the wavelet: in frequency, to build a nonparametric estimator of the spectral density from continuous time observations, and
in time, to bound the error of approximation of the wavelet coefficient with discrete time observations.
Conditions required on wavelet mothers are mild and satisfied by a large set of wavelets (Daubechies  wavelet $D_p$ for $p\ge 6$, Lemari\'e-Meyer, Morlet, Gabor, biorthogornal wavelets, \dots) and only exclude Haar basis and Daubechies  wavelet for $p\le 4$.
Actually, we do not need that the family of functions generated by dilations and translations forms a basis of $L^2(\R)$.
Secondly, our wavelet based estimator can be applied to stationary processes as well as processes having stationary increments. Moreover it is robust to eventual polynomial trends. Such properties are induced by the number of vanished moments of the mother wavelet. A periodogram estimator does not satisfy such  conditions and therefore can not be efficiently applied in so many cases.
\section{Main results}\label{maths}
This section contains three main results. In the first subsection, we specify conditions on the mother wavelet, and give a representation formula for the wavelet coefficients of the process. In the second subsection, we  establish a CLT satisfied by the sample variance wavelet coefficients. This result provides the rate of convergence of a spectral density estimator in parametric or semi-parametric cases (for instance for a fBm).
Eventually, the third subsection is devoted to a nonparametric estimation of the spectral density through a localization procedure.
\subsection{Definition and harmonizable representation of wavelet coefficients}
Let $\psi:\R \to \C$ be a function, the so-called "mother" wavelet, and denote $\ds\widehat{f}(\xi)=\int_{\R}
e^{-i\xi\, x}\,f(x)\,dx$ the Fourier transform of
$\ds f\in L^{1}(\R)\cap L^{2}(\R)$. Let $(m, q, r)\in \N^*\times \R_+^2$ and consider following set of
assumptions on $\psi$:
\vspace{0.3cm}
\newline {\bf Assumption W$(m, q, r)$}
 {\em $\psi:~\R \mapsto \C$ is a
differentiable function satisfying:
\begin{itemize}
\item {Number of vanishing  moments:}
for all $n\le m+1$, $\displaystyle{\int_{\R} \left |t^n\psi(t)\right |dt
<\infty}$, and
$$  \int_{\R}  t^n\psi(t)\,dt =0\quad for\;all\; n \le m.\quad$$
\item{Time localization:}
there exists a constant $C_{\psi} >0$ such that for all $t\in \R$,
\ban \big(1+|t|\big)^q \cdot
 \big| \psi(t)\big|
&\le& C_{\psi}.
 \ean
\item {Frequency localization:} there exists a constant $C'_{\psi} >0$ such that for all $\xi\in \R$,
\ban \big(1+|\xi|\big)^r \cdot
\big(\big|\widehat{\psi}(\xi)\big|+\big|\widehat{\psi}'(\xi)\big|\big)
&\le& C'_{\psi}.
 \ean
\end{itemize}}
\noindent
The first condition of W$(m, q, r)$ implies that
$\widehat{\psi}(\xi)=O(\xi^m)$ when $\xi\to 0$ and is
$(m+1)$ times continuously differentiable. In the sequel, we assume at most W($1,3,1/2$). These conditions
are mild and are satisfied by many famous wavelets (Daubechies  wavelet $D_p$ for $p\ge 6$, Lemari\'e-Meyer, Morlet, Gabor, biorthogornal wavelets, \dots).
It is also not mandatory to choose $\psi$ to be a ``mother''
wavelet associated to a multiresolution analysis of $\LL^2(\R)$ and the whole theory can be developed without resorting to this assumption.\\
~\\
Let $(a,b)\in \R_+^*\times \R$, and define $d_X(a,b)$ to be the
wavelet coefficient of the process $X$ for the scale $a$ and the
shift $b$, such that
$$
d_X(a,b):=\frac 1{\sqrt a} \int_{\R} \psi\left(
\frac{t-b}{a}\right)\,X(t)\,dt.
$$
This family of wavelet coefficients satisfies the following property:
\begin{Prop}[Harmonizable representation] \label{pte:d(a)}
Let $\psi$ satisfy Assumption W$(1, 1, 0)$ and $X$ be a Gaussian process
defined by (\ref{repr:harmonizable}) or (\ref{stat}) with a spectral density $f$
satisfying respectively (\ref{cond:f}) or (\ref{cond:f:stat}). Then,
\ba \label{repr:dX}&&
d_{X}(a, b) =\sqrt{a}\, \int_{\R} e^{ib\xi}\,
\overline{\widehat{\psi}}(a \xi)\, f^{1/2}(\xi)\,dW(\xi)~\mbox{for all $(a,b)\in \R_+^*\times \R$},
\ea
and, for $a>0$, $(d_{X}(a,b))_{b\in \R}$ is a stationary
centered Gaussian process with variance given by  \be
\label{def:I(a)} \E \left(\big |d_{X}(a,b) \big |^2\right)= \mathcal{I}_{1}(a) :=a \int_{\R} |\widehat{\psi}(au)|^2\,
f(u)\, du~\mbox{for all $b \in \R$}. \ee
\end{Prop}
The proof of this proposition is grouped with all the other proofs in the Appendix.
\subsection{An estimator of the variance of wavelet coefficients and its application to the semi-parametric estimation of the spectral density} \label{nonpar}
Let us begin with an example. If $X$ is a fBm with Hurst parameter $H\in (0,1)$, its spectral density is $f(\xi)=C \, |\xi|^{-(2H+1)}$ for all $\xi \in \R^*$ (with $C>0$). Then for a scale $a>0$ a straightforward computation of the variance of wavelet coefficients $\mathcal{I}_{1}(a)$ defined in (\ref{def:I(a)}) shows that $
\mathcal{I}_{1}(a)=K \, a^{2H+1}$ with $\ds K= C\, \Big (\int_{\R} \frac { |\widehat{\psi}(u)|^2}{u^{2H+1}}\,du\Big )$.
Therefore a consistent estimator of $\mathcal{I}_{1}(a)$ furnishes a consistent estimator of $H$ obtained by a log-log regression of $\ds \big (\mathcal{I}_{1}(a_i)\big)_{1\leq i\leq m}$ onto $(\log a_i)_{1\leq i \leq m}$. The same method works also for multiscale fBm (see Bardet and Bertrand, 2007b).
 \\
~\\
Thus our first aim is the estimation of $\mathcal{I}_{1}(a)$. When only  $\big(X(t^{(n)}_{0}),\ldots,X(t^{(n)}_{n})\big )$ is known, an explicit formula $\mathcal{I}_{1}(a)$ is not available for both the following reasons:
\begin{description}
\item [1.] on the one hand, $d_{X}(a,b)$ is defined with a Lebesgue integral and
cannot be directly computed from data. As in Gloter and Hoffmann (2007), an approximation formula will be considered for computing wavelet coefficients. Thus, for $(a,b)\in \R_+^*\times
\R$ we define the empirical wavelet coefficient by
\be
\label{def:eab}  e_{X}(a,b):=\frac {1}{\sqrt a} \sum
_{i=0}^{n-1}\Big (\int_{t^{(n)}_i}^{t^{(n)}_{i+1}}  \psi
\Big(\frac {t -b} a \Big)\, dt\Big ) \, X \big(t^{(n)}_i\big).
\ee
\item [2.] on the other hand, a sample mean of $|d_{X}(a,b) \big |^2$  instead on $\E \big |d_{X}(a,b) \big |^2$ is computable only. Thus, define the sample estimator of
$\mathcal{I}_{1}(a)$ by
\be \label{def:JN(a)} J_n(a) := \frac 1
{n+1} \,  \sum_{k=0}^{n}  \big |e_{X}(a,c_k)\big |^2, \ee where $(c_k)_{k}$ is a
family of increasing real numbers (so-called shifts).
In this paper, we will consider a uniform repartition of shifts, {\it i.e.} for $k=0,\ldots,n$,
\begin{eqnarray}\label{ck}
c_k=
T_n^\rho+k\, \frac {T_n-2\,T_n^\rho } n~~~\mbox{with $\rho \in (3/4,1)$}.
\end{eqnarray}
\end{description}
In this example $(c_k)_{1\leq k\leq n}$ are random variables depending on random times $(t^{(n)}_{1},\ldots,t^{(n)}_{n})$ but $c_{k+1}-c_k$ does not depend on $k$. We will see that it is not easy to consider a simpler expression of $(c_k)$; for instance $c_k=kT_n/n$ could not be used because there would be some edge effects for estimating the wavelet coefficients in $c_0$ or $c_n$. Therefore a sufficient ``distance'' from the boundaries  $0$ and $T_n$ is necessary. However, other choices of
$(c_k)_{k}$ are possible (for instance $c_k=t_k$) but we have not been able to find an optimal choice and
simulations do not show significant differences between these choices.
Now additional conditions on $f$ have to be considered:
\vspace{0.2cm}
\newline {\bf Assumption F($H$):} {\em $f$ is an even function, differentiable on $[0,\infty)$ except for a finite number $K$ of
real numbers $\omega_0=0<\omega_1<\dots<\omega_K$, but $f$ admits
left and right limits in $\omega_k$, with a derivative $f'$ (defined on all open
intervals $(\omega_{k},\omega_{k+1})$
with $\omega_{K+1}=\infty$ by convention) such that \be
\label{cond:f:prime} \int_{\R} \big(1\wedge |\xi|^3\big)\cdot
\left|f'(\xi)\right| \, d \xi < \infty. \ee Moreover, there exist
$C_0,C_0'>0$ and $H>0$, such that for all $|x|\ge \omega_K$
\ba \label{majoration:f:infini} f(x)\,\le\, C_0\,
|x|^{-(2H+1)}&\quad\mathrm{and}\quad& |f'(x)|\,\le\, C_0'\
|x|^{-(2H+2)}.\ea}
\noindent Here are several examples of processes having a spectral
density $f$ satisfying Assumption F($H$):\\
~\\
{\em Examples :} 1. A smooth Gaussian process having stationary increments satisfies F($H$) with $H\geq 1$ satisfies F($H$) with $H\ge1$.
\\
2. A fractional Brownian motion with Hurst parameter $H\in (0,\,1)$ satisfies F($H$). Indeed, its spectral density is given by
$f(\xi) =C\,  |\xi|^{-(2H+1)}\;$ and corresponds to a power law of the frequency.\\
3. However, fBm is a limited model. For instance, in some biological applications, 
statistical studies
suggest that the logarithm of the spectral density is a piecewise affine function of the log-frequency, see for instance Collins and De Luca (1993) or Billat {\em et al.} (2009). Furthermore in certain frequency bands the slope corresponds to a Hurst parameter $H$ larger than $1$.
 For these reasons, in Bardet and Bertrand (2007a), we have introduced the multiscale fBm such that there exists a family of frequencies $\omega_1<\dots<\omega_K$
satisfying $f(\xi) = C_i\,|\xi|^{-(2H_i+1)}$ for $|\xi| \in (\omega_i, \omega_{i+1})$ and $i=0,\dots,K$, with the convention that $\omega_0=0$ and $\omega_{K+1}=\infty$, $H_0<1$, $0<H_K$ and $(C_i,H_i)\in \R^*_+\times \R$ for $i=1,\ldots,K-1$.
Then Condition (\ref{cond:f}) and Assumption F($H$) are checked
with $H=H_K$.\\
4. A stationary process with a bounded spectral density  such as the Ornstein-Uhlenbeck process (for which $f(\xi):=\alpha \big ( \pi (\alpha^2+\xi^2)\big)^{-1}$ with $\alpha>0$). \\
~\\
The sample variance of wavelet coefficients $J_n(a)$ computed from the observed trajectory $(X (t^{(n)}_{0}),\ldots,X(t^{(n)}_{n}))$ and defined by (\ref{def:JN(a)}) satisfies the following CLT:
\begin{Theo} \label{Theo:TCL:discret}
Let $X$ be a Gaussian process defined by (\ref{repr:harmonizable}) or (\ref{stat})
with a spectral density $f$ satisfying (\ref{cond:f}) and Assumption F($H$),
$\psi$ satisfying Assumption W$(1,3,1)$ and $(c_k)_k$ defined by (\ref{ck}). If Assumption S($s$) holds with $\ds s> 2 +\frac 1 {2H} \, \big [1-3H\big]_+ $, and if
\ba \label{Cond:delta:T}
\E(T_n) \times \delta_n^{(s-1)\big (\frac {(2H\wedge 1)}{1+(2H \wedge 1)} \big ) \wedge \big ( \frac {1+(H\wedge 1)}{s+(H\wedge 1)}\big )} \limiten 0,
\ea
then for all $a>0$,
\begin{eqnarray} \label{TLC2d}
~~~~~~~~\sqrt{\E \, T_n} \,  \big ( J_n(a)-\mathcal{I}_1(a) \big )
\limiteloin \mathcal{N}\Big (0\, ,\, 4 \pi \,
a^2   \, \int_{\R}\big |
\widehat{\psi}(a z)\big |^4 f^2(z)\, d z \Big ).
\end{eqnarray}
\end{Theo}
\begin{rem}\begin{enumerate}
             \item The convergence rate of the CLT (\ref{TLC2d}) is $\sqrt{\E(T_n)}$ when Condition  (\ref{Cond:delta:T}) is satisfied. A natural question is what happens elsewhere? This leads to the following comments: roughly speaking, from theorem \ref{Theo:TCL:discret}  and lemma \ref{lem:maj:vn:stoch}, one can deduce
\ba
\label{TCL:heuristic}
J_n(a)&=&\mathcal{I}_1(a)+ \big[\E(T_n)\big]^{-1/2}\,\Gamma\, U + \zeta_n
\ea
where $U\sim \mathcal{N}(0,1)$, $\Gamma^2=4 \pi \,
a^2   \, \int_{\R}\big |
\widehat{\psi}(a z)\big |^4 f^2(z)\, d z$ corresponds to the variance in CLT (\ref{TLC2d}) and $\zeta_n$ corresponds to the discretization error.
As soon as Condition (\ref{Cond:delta:T}) is fulfilled, the discretization term $\zeta_n$ is negligible with respect to the CLT term 
and the rate of convergence is $\sqrt{\E(T_n)}$.
If the condition (\ref{Cond:delta:T}) is not satisfied, then the upper bound of the mean square error  does no longer decrease when $\E(T_n)\to \infty$.
\item It also possible to specify Condition (\ref{Cond:delta:T}) by
using a relation between $\delta_n$ and $n$; for this, let
$$
\delta_n=C_\delta \, n^{-d}\quad \mbox{with $0<d<1$}.
$$
The following Table 1  
summarizes the possible choices of $s$ and $d$  and the supremum of  the
convergence rate of the CLT (\ref{TLC2d})) for several cases.
\end{enumerate}
\end{rem}

\begin{table}[h]\label{Table0}
{
\begin{center}
\small
\begin{tabular}{|c|c|c|c|c|}
&$H$ known &  $H$ unknown &  $H$ unknown  & $H$ unknown   \\
& $H>0$ &  $H>0$& $H\geq 1/3$ & $H\geq 1$ \\
\hline \hline
Condition on $s$  &$s>  2 +\frac 1 {2H} \, \big [1-3H\big]_+ $ & $s=\infty$& $s >2$ & $s>2$ \\
\hline Condition on $d$ & $d >\big ( \frac {1+(2H\wedge 1)}{1+s(2H\wedge 1)}\big )\vee \big (\frac {s+(H\wedge 1)}{s(2+(H\wedge 1))-1}  \big )$&$d \geq  \frac 1 2 $ &$ d>(\frac 2 {1+s})\vee (\frac{2s+1}{5s-2}) $ &$d>(\frac 2 {1+s})\vee (\frac{s+1}{3s-1}) $ \\
\hline  Supremum of the rate of & $n^{\frac{1-d}2 } $ for $s>  2 +\frac 1 {2H} \, \big [1-3H\big]_+ $ &$ n^{-\frac 1 4 }$ & $n^{-\frac 1 4 }$ for $s=4$ &$ n^{-\frac 1 4 }$ for $s=3$\\
convergence of CLT (\ref{TLC2d})&$n^{-\frac 1 2 \frac {1+ (H\wedge 1)}{2+(H\wedge 1)} } $  for $s=\infty$ &  & $n^{-\frac 3 {10} }$ for $s=\infty$ &$ n^{-\frac 1 3 }$ for $s=\infty$ \\
\hline  Supremum of the rate of & $n^{2\frac{1-d}5 } $ for $s>  2 +\frac 1 {2H} \, \big [1-3H\big]_+ $ &$ n^{-\frac 1 5 }$ & $n^{-\frac 1 5 }$ for $s=4$ &$ n^{-\frac 1 5 }$ for $s=3$\\
convergence of CLT (\ref{Jpp})&$n^{-\frac 2 5 \frac {1+ (H\wedge 1)}{2+(H\wedge 1)} } $  for $s=\infty$ &  & $n^{-\frac 6 {25} }$ for $s=\infty$ &$ n^{-\frac 4 {15} }$ for $s=\infty$\\
\hline
\end{tabular}
\end{center}}
\caption {\it Conditions on $s$, $d$ and supremum of the convergence rate of the semiparametric estimator of the spectral density $f$ (CLT (\ref{TLC2d})) and the nonparametric estimator of $f$ (CLT (\ref{Jpp}) with $\lambda_n \simeq C T_n^{1/5+\kappa}\simeq C n^{1/5(1-d)+\kappa}$ with $\kappa>0$ arbitrarily small) following the a priori on $H$.}
\end{table}
Note that   CLT (\ref{TLC2d}) can be applied to an estimation of each $H_i$ of a multiscale fractional Brownian motion when a trajectory is observed at random times. Indeed, in such a case and if $\psi$ is chosen such that $\widehat \psi (\xi)\neq 0$ only for $\xi \in [-\beta,-\alpha]\cup [\alpha,\beta]$, then (see details in Bardet and Bertrand, 2007b):
$$
\mathcal{I}_1(a)=a^{2H_i+1} C_{f,\psi}~~~\mbox{for all $\xi \in [\alpha/\omega_i,\beta/\omega_{i+1}]$},
$$
with $C_{f,\psi}>0$ not depending on $a$. Therefore a log-log-regression of $  J_n(a)$ onto $a$ for several values of $a \in [\alpha/\omega_i,\beta/\omega_{i+1}]$ provides an estimator of $H_i$ and $C_i$ which follows a CLT with the same convergence rate as (\ref{TLC2d}). Such a result may of course also be applied to fBm without specifications of the scales $a$. This is more precisely stated in the following
 \begin{cor}[parametric case]\label{cor_loc}
Let $X$ be a  fBm with parameters $H\in (0,1)$ and $C>0$. Assume that $(X_{t_1^{(n)}},\cdots,  X_{t_1^{(n)}})$ is observed, that Assumption S($\infty$) is fulfilled, that $\psi$ satisfies Assumption $W(1, 3, 1)$ and that $\delta_n \E T_n \limiten 0$.
 Then there exists a constant $c>0$ such as for $n$ large enough,
\ba\label{estim_para}
\E \Big [\Big \|\Big ( \begin{array}{c} \widehat H_n \\ \widehat C_n \end{array} \Big )-\Big ( \begin{array}{c} H \\ C \end{array} \Big ) \Big \|^2\Big ] \leq \frac c {\E T_n}.
\ea
where $\widehat H_n$ and $\widehat C_n$ are the estimators obtained by log-log-regression of $  J_n(a)$ onto $a$.
If moreover, the Hurst index $H$ is known in advance to lie in the interval $(1/3, 1)$, then Condition  S($\infty$) can be replaced by Condition  S($2+\varepsilon$) for any $\varepsilon>0$.
\end{cor}
 To our knowledge, only Begyn (2005) provides an asymptotic result on the estimation of $H$ under irregular observation times but only in the case of fBm and with a stronger condition than Assumption S($\infty$).
\subsection{A nonparametric estimator of the spectral density}
The third result of this paper deals with the pointwise estimation of $f$ through a localization procedure in theorem \ref{Theo:TCL:discret}. Let us define the "rescaled" functions:
\ba
\label{def:recale:psi}
\psi_\lambda(x) &:=&  \frac 1 {\sqrt{\lambda}}\, e^{i\,x}\,
\psi\big ( \frac x \lambda \big )
\ea
in a way that enables the convergence, as $\lambda\to \infty$, of $|\widehat{\psi}_\lambda|^2$ to a Dirac mass concentrated at the frequency $\xi=1$. Then a rescaled version of the estimator (\ref{def:eab}, \ref{def:JN(a)}) is introduced:
\ba
\label{def:rescaled:estimator}
\widehat f_n^{(\lambda)}(\xi)
&:=& \frac \xi {\left\|\psi\right\|^2_{{\cal L}^2}} \,
\frac {1}{n+1} \,  \sum_{k=0}^{n}
\Big |\sum_{i=0}^{n-1} X \big(t^{(n)}_i\big) \, \int_{t^{(n)}_i}^{t^{(n)}_{i+1}}  \psi_{\lambda}
\big(\xi \big (t -c_k\big ) \big)\, dt \Big | ^2.
\ea
From (\ref{def:recale:psi}), it is obvious that
$$
\widehat{\psi}_\lambda(\xi) = \sqrt{\lambda}\,
\widehat{\psi}\big(\lambda(\xi-1)\big)~~~~\mbox{$\forall \xi \in \R$},
$$
and after that $$\ds{\cal I}_\lambda(a):=\int_{\R} |\widehat{\psi}_\lambda(u)|^2\,
f(u/a)\, du \to f(1/a)\, \left\|\psi\right\|^2_{{\cal L}^2}\quad \mathrm{when} \quad\lambda\to \infty$$ under weak conditions.
Then a CLT is established for the nonparametric estimator (\ref{def:rescaled:estimator}) with a sequence $(\lambda_n)_n$ satisfying  $\lambda_n\to \infty$ and under the assumptions of theorem \ref{Theo:TCL:discret}.
Note that the first condition $\psi_\lambda  \in W(1,3,1)$ is fulfilled as soon as $\lambda_n>\Lambda$ when $\widehat \psi$ is compactly supported in $[-\Lambda, \Lambda]$. Now, by using an appropriate choice
of a sequence $(\psi_{\lambda_n})$, one obtains:
\begin{Prop}\label{cvg:Ip}
Assume that the assumptions of theorem \ref{Theo:TCL:discret} hold. If the spectral density $f$ is a twice continuously differentiable function on $\R^*$, if $\widehat \psi$ is compactly supported, and if the sequence  $(\lambda_n)_n$ is such that $\ds\frac {\lambda_n^2}{n\delta_n} \limiten 0$ and $\ds\frac {\lambda_n^5}{n\delta_n} \limiten \infty$, then for all $ \xi>0$,
\ba \label{Jpp}
~~~~\sqrt {\frac{T_n} {\lambda_n} }\, \big (\widehat f_n^{(\lambda_n)}(\xi)
- f(\xi) \big ) \limiteloin \mathcal{N}\Big (0\, ,\, \frac {4 \, \pi} \xi \,f^2(\xi) \,\frac { \int_{\R}\big |
\widehat{\psi}(u)\big |^4 \, d u} {\big ( \int_{\R}\big |
\widehat{\psi}(u)\big |^2 \, d u\big )^2} \Big ).
\ea
\end{Prop}
The rate of convergence of the parametric (or semiparametric) estimator is $T_n^{-1/2}$, see  CLT (\ref{TLC2d}). In the case of a nonparametric estimator, using the optimal choice of $\lambda_n$, {\it i.e.} $\lambda_n=C(n\delta_n)^{1/5+\kappa}=O(T_n^{1/5+\kappa})$ with $\kappa>0$ arbitrary small, the supremum of the convergence rate of the estimator is $T_n^{-2/5}$. This is the same rate of convergence as for the periodogram of a stationary process observed in continuous time (Parzen, 1983) or observed during random times satisfying Masry's conditions (Lii and Masry, 1994). However, in this last case, $T_n \sim C\, n$ p.s. when $n \to \infty$. Our result, that is CLT (\ref{Jpp}), is clearly more general: it concerns processes having stationary increments and satisfying weak conditions on the random observation times. But the prize to pay for obtaining the convergence rate $T_n ^{-2/5}$ is that $T_n \sim  C\, n^{1-d}$ with $\ds d >\big ( \frac {1+(2H\wedge 1)}{1+s(2H\wedge 1)}\big )\vee \big (\frac {s+(H\wedge 1)}{s(2+(H\wedge 1))-1}  \big )$, {\it i.e.} for instance $T_n \sim  C\, n^{\frac 1
2}$ for $s=3$ and $H\geq 1$, or for $s=\infty$ and $H>0$; more details are provided in Table~1.

How to explain that this convergence rate depends on $s$ and $H$? On the one hand, the smaller $H$ the more irregular the trajectory of $X$ when $X$ is a process having stationary increments (the Hölder parameter of a trajectory of $X$ is then $H^+$ for all $H^+<H$). Therefore empirical wavelet coefficients, defined almost as Riemann sums approximate better a smooth path than an irregular path and this explains that the smaller $H$ the smaller the convergence rate of CLT (\ref{Jpp}). In the case of stationary processes, $H\geq 1$ and the convergence rate does not depend on $H$.
On the other hand, the smaller $s$ the more variable the observed times. Then for each frequency there are not enough successive data with appropriate lag allowing to correctly estimate the spectral density around this frequency. Then the smaller $s$ the smaller the convergence rate of CLT (\ref{Jpp}).

Moreover, under W($m,3, 1$), $\ds\int t^n\psi(t)dt=0$ for all $n \leq m$ and all wavelet coefficients of any polynomial function with degree less than or equal to $m$ vanish. Therefore, the estimator $\widehat f_n^{(\lambda_n)}$ is robust, since:
\begin{cor}\label{corpoly}
Under Assumption W($m, 3, 1$) with $m \in \N^*$, Proposition \ref{cvg:Ip} holds when a polynomial trend with degree less than or equal to $m$ is added to $X$.
\end{cor}
\section{Numerical experiments}\label{numeric}
For the  numerical applications, one has chosen:
\begin{enumerate}
\item $\psi$ is chosen such that $\widehat{\psi}(\xi)=\exp \big (-(|\xi|\cdot(5-|\xi|))^{-1}\big )\1_{|\xi|\leq 5}(\xi)$, which satisfies
Assumption W($m,r$) for any $(m,r)$ (and $\widehat \psi(\xi)=0$ for $|\xi|\geq 5$).
\item $\delta_n=n^{-0.6}$ for insuring the convergence of $\widehat f_n^{(\lambda_n)}(\xi)$ for any $H>0$ and $s\geq 3$.
\item $\lambda_n=n^{d'}$ with $1/6 <d'<1/2$. The admissibility condition on wavelets $(\psi_{\lambda_n})$ requires that $n^{d'}\geq \Lambda=5$ and after numerous simulations, we have chosen $d'=\log(15)/\log(n)$.
\end{enumerate}
\subsection{Estimation of the spectral density of a fractional Brownian motion observed at random times}
For a standard ($\E X^2(1)=1$) fBm with Hurst parameter $H$, $f(\xi)=C(H)\, |\xi|^{-2H-1} d\xi $
with $C(H)=\big (H\Gamma(2H)\sin(\pi H) \big)/\pi$.
Three different kinds of independent and identically distributed random times are considered:
\begin{description}
\item {\bf (T1)}: non-random uniform sampling, such that $L_k=1$ for all $k \in \N^*$;
\item {\bf (T2)}: exponential random times, such that $\E L_k=1$ for all $k \in \N^*$;
\item {\bf (T3)}: random times such that for for all $k \in \N^*$, the cumulative distribution
function of $L_k$ is $F_{L_k}(x)=(1-x^{-4})\1_{x\geq 1}$ implying $\E L_k^p<\infty$ for all $p<4$ and $\E L_k^4=\infty$. In this case Assumption S($s$) is satisfied if and only if $s<4$.
\item {\bf (T4)}: random times such that for for all $k \in \N^*$, the cumulative distribution
 function of $L_k$ is $F_{L_k}(x)=(1-x^{-2})\1_{x\geq 1}$ implying $\E L_k^p<\infty$ for all $p<2$ and $\E L_k^2=\infty$. In this case Assumption S($s$) is satisfied if and only if $s<2$. As a consequence, this case does not satisfy the hypothesis of theorem \ref{Theo:TCL:discret}.
\end{description}
An example of such estimation of the spectral density for $H=0.2$, $N=50000$ and random times T2 is presented in Figure 1. The results of the simulations are also provided in Table 2. \\
~\\
{\bf Fig.~1} {\it An example of the estimation of the spectral density (left) and its logarithm (right) of a FBM observed at exponential random times (T2) with confidence intervals ($H=0.2$, $N=50000$).} 
\[
\epsfxsize 7cm \epsfysize 6cm \epsfbox{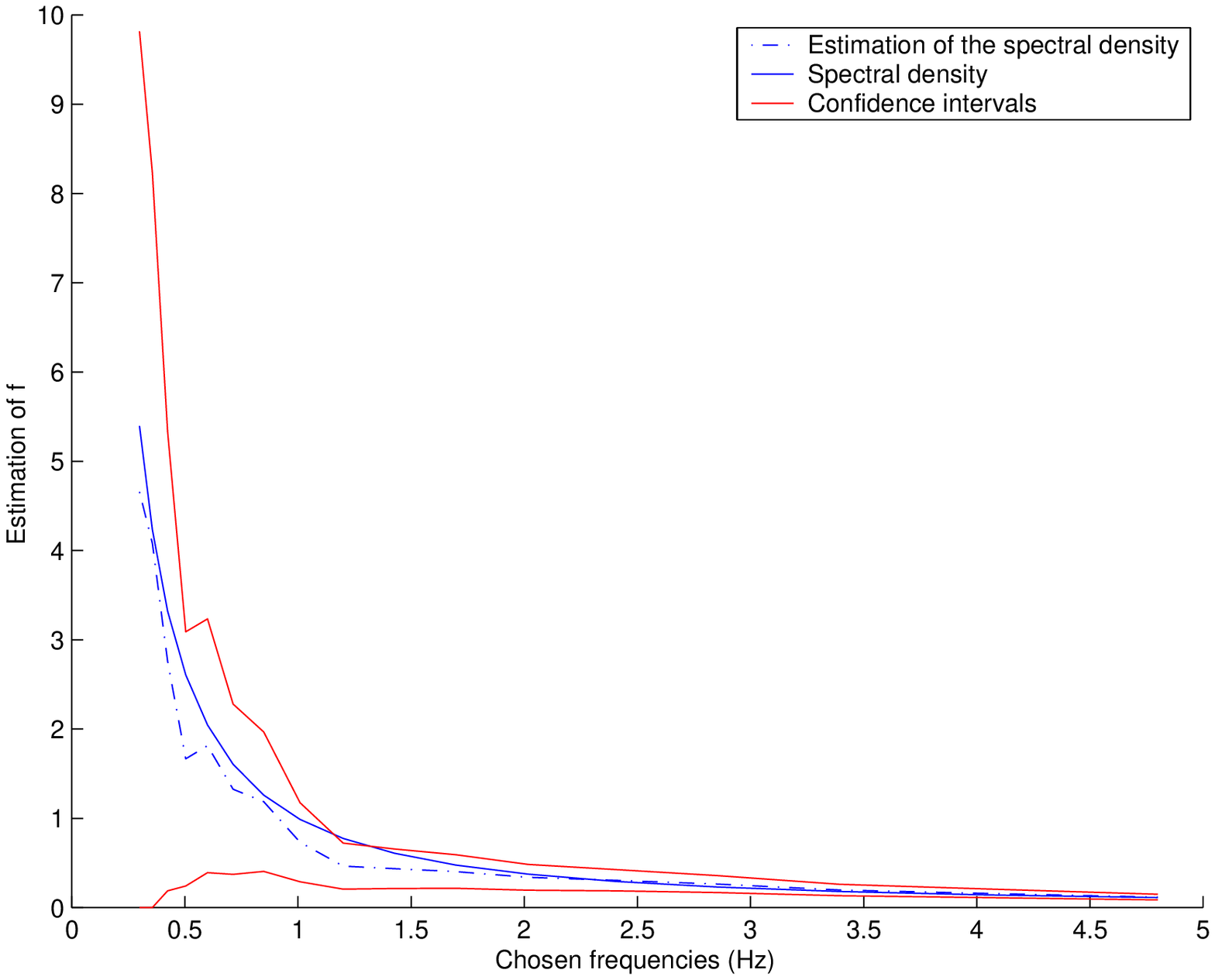}
\hspace*{1.5 cm} \epsfxsize 7cm \epsfysize 6cm
\epsfbox{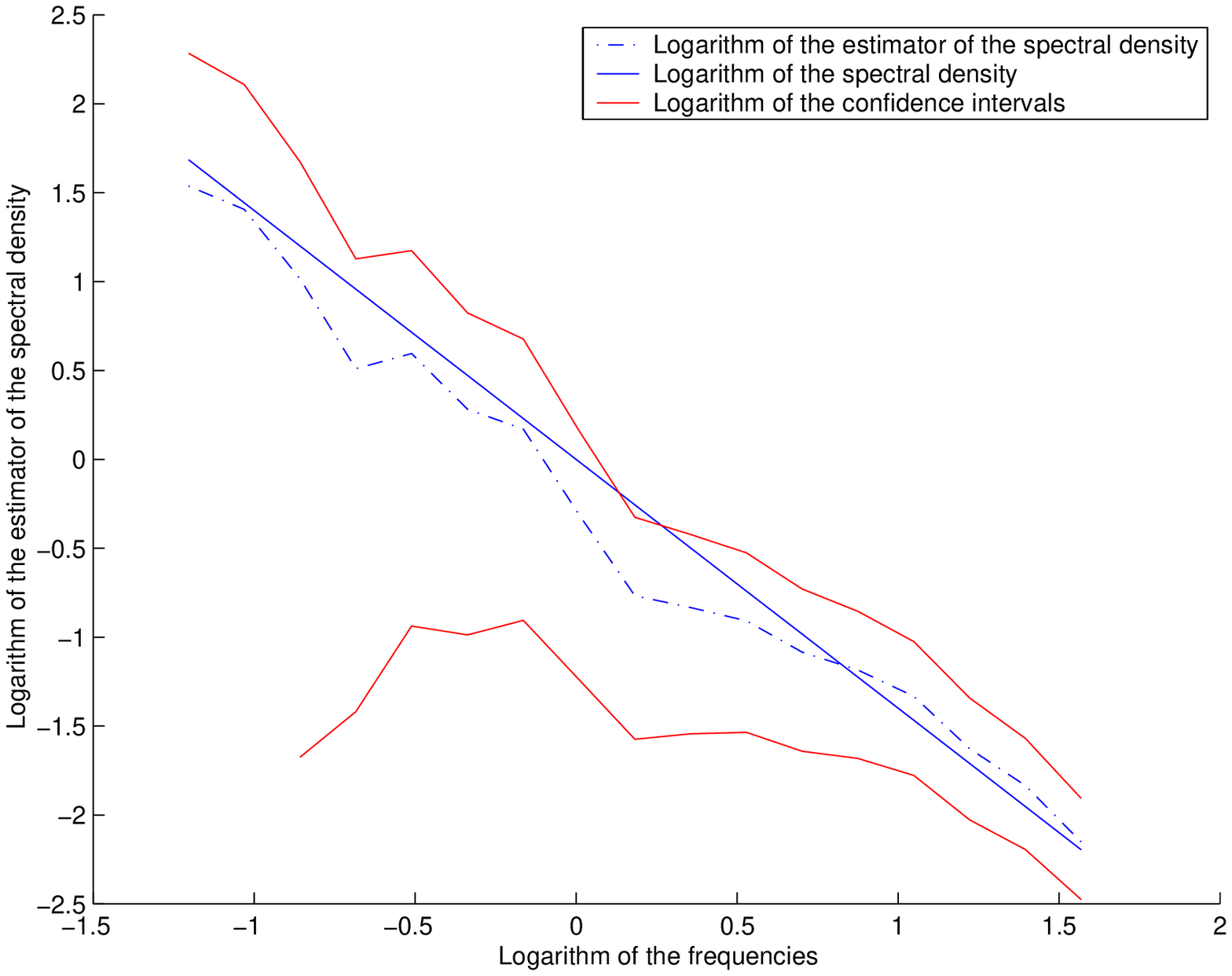}
\]

\begin{table}[h]\label{Table2}
{
\begin{center}
$N=10^3 $ ~\begin{tabular}{|c|c|c|c|c|}
 &  & $H=0.2$ & $H=0.5$ & $H=0.8$  \\
\hline \hline
T1& $\sqrt{MSE}$ of ${\widehat f}_N(1)$ & 0.47&0.65&0.77 \\
 &  $\widehat{MISE}$ on $[0.3,5]$ & 2.53 & 13.50 & 80.89\\
\hline
T2& $\sqrt{MSE}$ of ${\widehat f}_N(1)$ & 0.65&0.67&0.75  \\
 &  $\widehat{MISE}$ on $[0.3,5]$ & 3.64 & 10.65 & 39.85\\
\hline
T3& $\sqrt{MSE}$ of ${\widehat f}_N(1)$ & 0.42&0.72&1.20  \\
 &  $\widehat{MISE}$ on $[0.3,5]$ & 2.48 & 7.83 & 55.20\\
\hline
T4& $\sqrt{MSE}$ of ${\widehat f}_N(1)$ & 1.03&3.34&2.44 \\
&  $\widehat{MISE}$ on $[0.3,5]$ & 6.07 & 84.05 & 144.40\\
\hline
\hline
\end{tabular}
\end{center}
\begin{center}
$N=10^4 $ ~\begin{tabular}{|c|c|c|c|c|}
 &  & $H=0.2$ & $H=0.5$ & $H=0.8$  \\
\hline \hline
T1& $\sqrt{MSE}$ of ${\widehat f}_N(1)$ & 0.35&0.37&0.79  \\
 &  $\widehat{MISE}$ on $[0.3,5]$ & 0.95 & 3.90 & 57.19\\
\hline
T2& $\sqrt{MSE}$ of ${\widehat f}_N(1)$ & 0.45&0.47&0.29  \\
 &  $\widehat{MISE}$ on $[0.3,5]$ & 1.04 & 3.17 & 16.26\\
\hline
T3& $\sqrt{MSE}$ of ${\widehat f}_N(1)$ & 0.47&0.46&0.95  \\
 &  $\widehat{MISE}$ on $[0.3,5]$ & 1.20 & 4.91 & 26.6\\
\hline
T4& $\sqrt{MSE}$ of ${\widehat f}_N(1)$ & 0.61&0.61&1.74 \\
&  $\widehat{MISE}$ on $[0.3,5]$ & 2.74 &9.55  & 49.55\\
\hline
\hline
\end{tabular}
\end{center}
\begin{center}
\hspace{-0.4cm} $N=5 \cdot 10^4 $ ~\begin{tabular}{|c|c|c|c|c|}
 &  & $H=0.2$ & $H=0.5$ & $H=0.8$  \\
\hline \hline
T1& $\sqrt{MSE}$ of ${\widehat f}_N(1)$ & 0.36&0.30&0.40 \\
 &  $\widehat{MISE}$ on $[0.3,5]$ & 0.81 & 2.60& 10.77\\
\hline
T2& $\sqrt{MSE}$ of ${\widehat f}_N(1)$ & 0.21&0.22&0.31  \\
 &  $\widehat{MISE}$ on $[0.3,5]$ & 1.07 & 2.07 & 7.65\\
\hline
T3& $\sqrt{MSE}$ of ${\widehat f}_N(1)$ & 0.34&0.26&0.48  \\
 &  $\widehat{MISE}$ on $[0.3,5]$ & 0.74 & 3.17 & 13.3\\
\hline
T4& $\sqrt{MSE}$ of ${\widehat f}_N(1)$ & 0.40&0.56&2.59 \\
&  $\widehat{MISE}$ on $[0.3,5]$ & 1.02 &5.69  & 41.41\\
\hline
\hline
\end{tabular}
\end{center}
}
\caption{\it Consistency of the estimator ${\widehat f}_N$ in the case of paths of a FBM observed at random times ($50$ independent replications are generated in each case).}
\end{table}
~\\
\noindent Comments on simulation results: 
\begin{enumerate} 
\item The larger $N$ the more accurate the estimator of $f$
except for the case of random times T4 (a case not included in the conditions of Proposition \ref{cvg:Ip});
\item The results are similar for T1 and T2, a little less accurate for T3; 
\item The smaller $H$ the more accurate the estimator of $f$.
\end{enumerate}
\subsection{Estimation of the spectral density of a stationary Ornstein-Uhlenbeck process}
Here, instead of FBM which is a process having stationary increments, we consider a stationary Ornstein-Uhlenbeck process which is a Gaussian stationary process with covariance $r(t):=\exp(-\alpha |t|)$ and therefore with spectral density $f(\xi):=\alpha \big ( \pi (\alpha^2+\xi^2)\big)^{-1}$. In such a case, since the spectral density is an analytic function, there exists a more accurate nonparametric estimator (see for instance, Ibragimov, 2004). However, to our knowledge, the case when paths are observed at random times is not considered in the literature. The results of simulations are provided in Table 3. \\
\begin{table}[t]\label{Table3}
{
\small
\begin{center}
$N=10^3 $ ~\begin{tabular}{|c|c|c|c|c|}
 &  & $\alpha=0.1$ & $\alpha=1$ & $\alpha=10$  \\
\hline \hline
T1& $\sqrt{MSE}$ of ${\widehat f}_N(0.3)$ & 0.51&0.22&0.020 \\
 &  $\widehat{MISE}$ on $[0.3,5]$ & 0.022 & 0.014 & 0.00067\\
\hline
T2& $\sqrt{MSE}$ of ${\widehat f}_N(0.3)$ & 0.30&0.30&0.021  \\
 &  $\widehat{MISE}$ on $[0.3,5]$ & 0.010 & 0.024 & 0.0010\\
\hline
T3& $\sqrt{MSE}$ of ${\widehat f}_N(0.3)$ & 0.36&0.23&0.018 \\
 &  $\widehat{MISE}$ on $[0.3,5]$ & 0.00052 & 0.015 & 0.00052\\
\hline
T4& $\sqrt{MSE}$ of ${\widehat f}_N(0.3)$ & 0.28&0.23&0.032 \\
&  $\widehat{MISE}$ on $[0.3,5]$ & 0.016 & 0.016 & 0.0045\\
\hline
\hline
\end{tabular}
\end{center}
\begin{center}
$N=10^4 $ ~\begin{tabular}{|c|c|c|c|c|}
&  & $\alpha=0.1$ & $\alpha=1$ & $\alpha=10$  \\
\hline \hline
T1& $\sqrt{MSE}$ of ${\widehat f}_N(0.3)$ & 0.20&0.18&0.017  \\
 &  $\widehat{MISE}$ on $[0.3,5]$ & 0.0033 & 0.0088& 0.00031\\
\hline
T2& $\sqrt{MSE}$ of ${\widehat f}_N(0.3)$ & 0.14&0.18&0.019  \\
 &  $\widehat{MISE}$ on $[0.3,5]$ & 0.0032 & 0.0092 & 0.00036\\
\hline
T3& $\sqrt{MSE}$ of ${\widehat f}_N(0.3)$ & 0.17&0.18&0.016 \\
 &  $\widehat{MISE}$ on $[0.3,5]$ & 0.0027 & 0.011 & 0.00032\\
\hline
T4& $\sqrt{MSE}$ of ${\widehat f}_N(0.3)$ & 0.18&0.13&0.024\\
&  $\widehat{MISE}$ on $[0.3,5]$ & 0.0058 &0.0095  & 0.00037\\
\hline
\hline
\end{tabular}
\end{center}
\begin{center}
\hspace{-0.4cm} $N=5 \cdot 10^4 $ ~\begin{tabular}{|c|c|c|c|c|}
&  & $\alpha=0.1$ & $\alpha=1$ & $\alpha=10$  \\
\hline \hline
T1& $\sqrt{MSE}$ of ${\widehat f}_N(0.3)$ & 0.14&0.10&0.012 \\
 &  $\widehat{MISE}$ on $[0.3,5]$ & 0.0016 & 0.0045& 0.00015\\
\hline
T2& $\sqrt{MSE}$ of ${\widehat f}_N(0.3)$ & 0.26&0.13&0.011 \\
 &  $\widehat{MISE}$ on $[0.3,5]$ & 0.012 & 0.0055 & 0.00014\\
\hline
T3& $\sqrt{MSE}$ of ${\widehat f}_N(0.3)$ & 0.18&0.14&0.012  \\
 &  $\widehat{MISE}$ on $[0.3,5]$ & 0.0023 & 0.0049 & 0.00017\\
\hline
T4& $\sqrt{MSE}$ of ${\widehat f}_N(0.3)$ & 0.16&0.16&0.017 \\
&  $\widehat{MISE}$ on $[0.3,5]$ & 0.0084 &0.034 & 0.00019\\
\hline
\hline
\end{tabular}
\end{center}
}
\caption{\it Consistency of  ${\widehat f}_N$ in the case of paths of a stationary Ornstein-Uhlenbeck process observed at random times ($50$ independent replications are generated in each case).}
\end{table}
\newline \noindent Comments on simulation results: 
\begin{enumerate}
\item The larger $N$ the more accurate the estimator of $f$ for all choice of random time
\item The results are similar for T1, T2, T3
and a less accurate for T4;
\item Surprisingly, the case $\alpha=1$ is not clearly better than $\alpha=0.1$ despite the fact that the larger $\alpha$ the less correlated the process.
\end{enumerate}
\subsection{Estimation of the spectral density of heartbeat time series}\label{cardio}
Heartbeats of several working people have been recorded during $24$ hours (see an example in Fig. 2). These data have kindly been furnished by professors Alain Chamoux and Gil Boudet (Faculty of Medicine, Occupational safety and health,  University of Auvergne, Clermont-Ferrand). \\
~\\
{\bf Fig. 2} {\it An example of heart inter-beats during $24h$}
\[
\epsfxsize 12cm \epsfysize 5cm \epsfbox{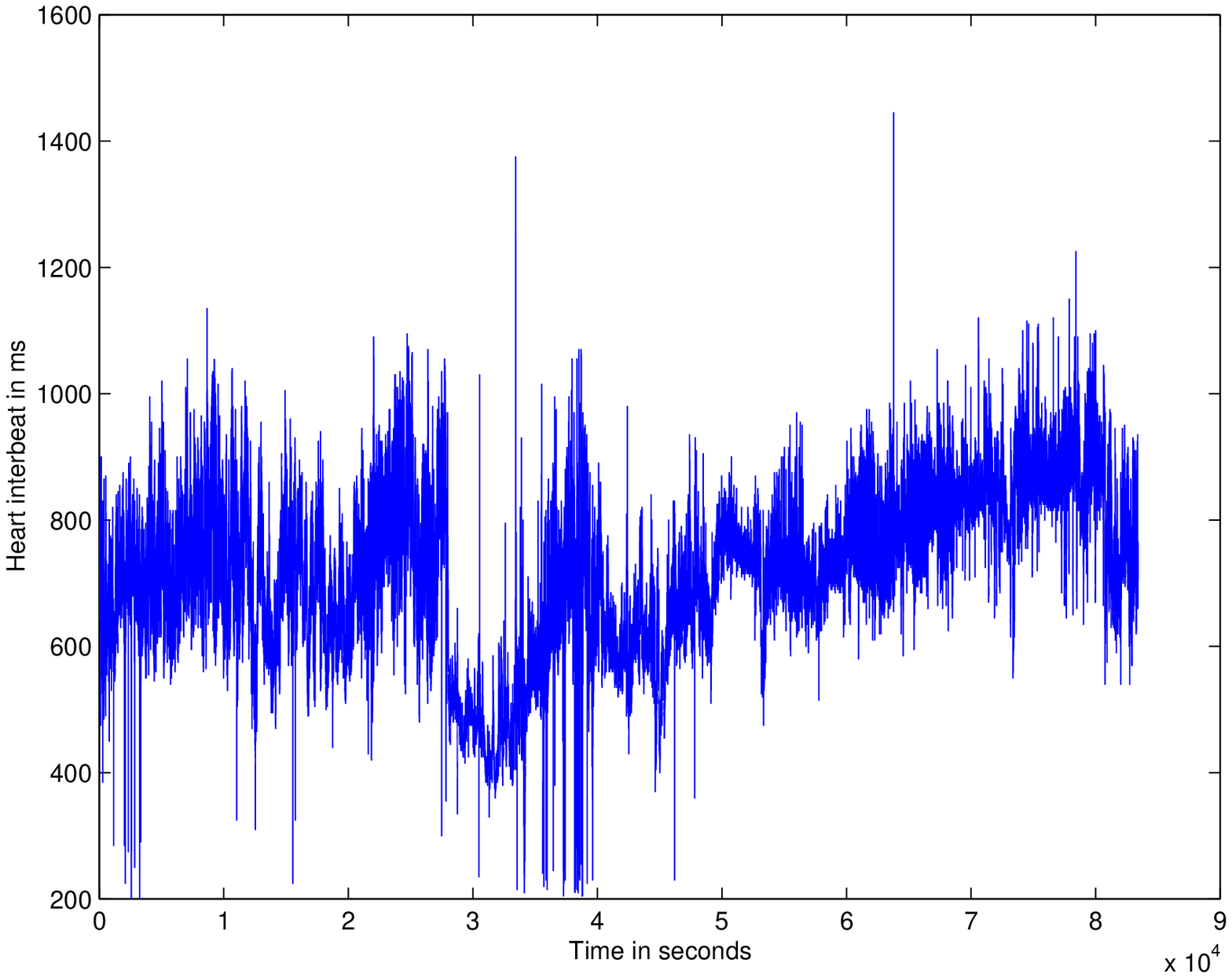}\]
Actually, heartbeat is measured at frequency $100$Hz. 
For physiological reasons, the duration between two observations should be between $250$ milliseconds and $2$ seconds. In this framework, Assumption S($\infty$) would  obviously be satisfied.
Assumption S($s$) with $s<\infty$ holds for instance for other physiological signals like EMG, EEG\dots

Cardiologists are interested in the study of this signal in two frequency
bands: the orthosympathetic and parasympathetic bands, i.e., the frequency bands
$(0.04\,Hz,\,0.15\, Hz)$ and $(0.15\, Hz,\,0.5\, Hz)$ respectively.
The definition of these bands is the outcome of research works, see e.g., Task
force of the European Society of Cardiology and the North American Society of
Pacing and Electrophysiology (1996), and is based on
the fact that the behavior of the energy contained inside these bands would be a relevant
indicator on the level of the stress of an individual.

Indeed, for the heart rate, the parasympathetic system is often compared to
the brake while the orthosympathetic system would be a nice accelerator; see
e.g. Goldberger (2001). At rest
there is a permanent braking effect on the heart rate. Any solicitation of the
cardiovascular system, any activity initially produces a reduction of
parasympathetic brake followed by a gradual involvement of the sympathetic
system. These mechanisms are very interesting to watch in many diseases
including heart failure, but also rhythm disorders that may fall under one or
other of these two effects, monitoring the therapeutic effect of several
medicines including some psychotropic. In the field of physiology such data
are crucial for measuring the level of stress induced by
physical activity or level of perceived stress, which can be considered as a
criterion of overtraining in sport.

We decompose these data in $3$ temporal zones following the activity:
\begin{itemize}
\item Quiet activities ($t \in [1,28000]$ in seconds);
\item Intensive activities ($t \in [28000,51400]$ in seconds);
\item Sleep ($t \in [60000,83400]$ in seconds).
\end{itemize}
Applying the spectral density estimator on those $3$ sub-data sets and plotting its log-log representation for frequencies in $[0.02,1]$ Hz, we observe that:
\begin{itemize}
\item in zone ``Sleep'' (see Figure 3), only one regression line could be computed for frequencies in $[0.04,0.5]$ Hz which is the usual spectral interval considered by specialists; in this zone $\widehat H \simeq 0.99$;
\item in zone ``Quiet activities'' (respectively ``Intensive activities''), (see Figure 3), two regression lines could be drawn for frequencies in $[0.04,0.5]$ Hz, distinguishing the orthosympathetic and the parasympathetic spectral domains. Using an algorithm computing the ``best'' two regression lines (see for instance Bardet and Bertrand, 2007b), one obtains that $H\simeq 1.34$ (respectively $H \simeq 1.44$) in the orthosympathetic domain which is $[0.04,0.09]$ Hz (respectively $[0.04,0.11]$ Hz) and $H\simeq 0.89$ (respectively $H \simeq 0.79$) in the parasympathetic domain which is $[0.09,0.5]$ Hz (respectively $[0.11,0.5]$ Hz).
\end{itemize}
{\bf Fig. 3} {\it Log-log representation of the spectral density estimator during ``Sleep'' zone (upper), ``Quiet'' activities (middle) and ``Intensive activities'' zone (lower)}
\[
\epsfxsize 12cm \epsfysize 6.5cm \epsfbox{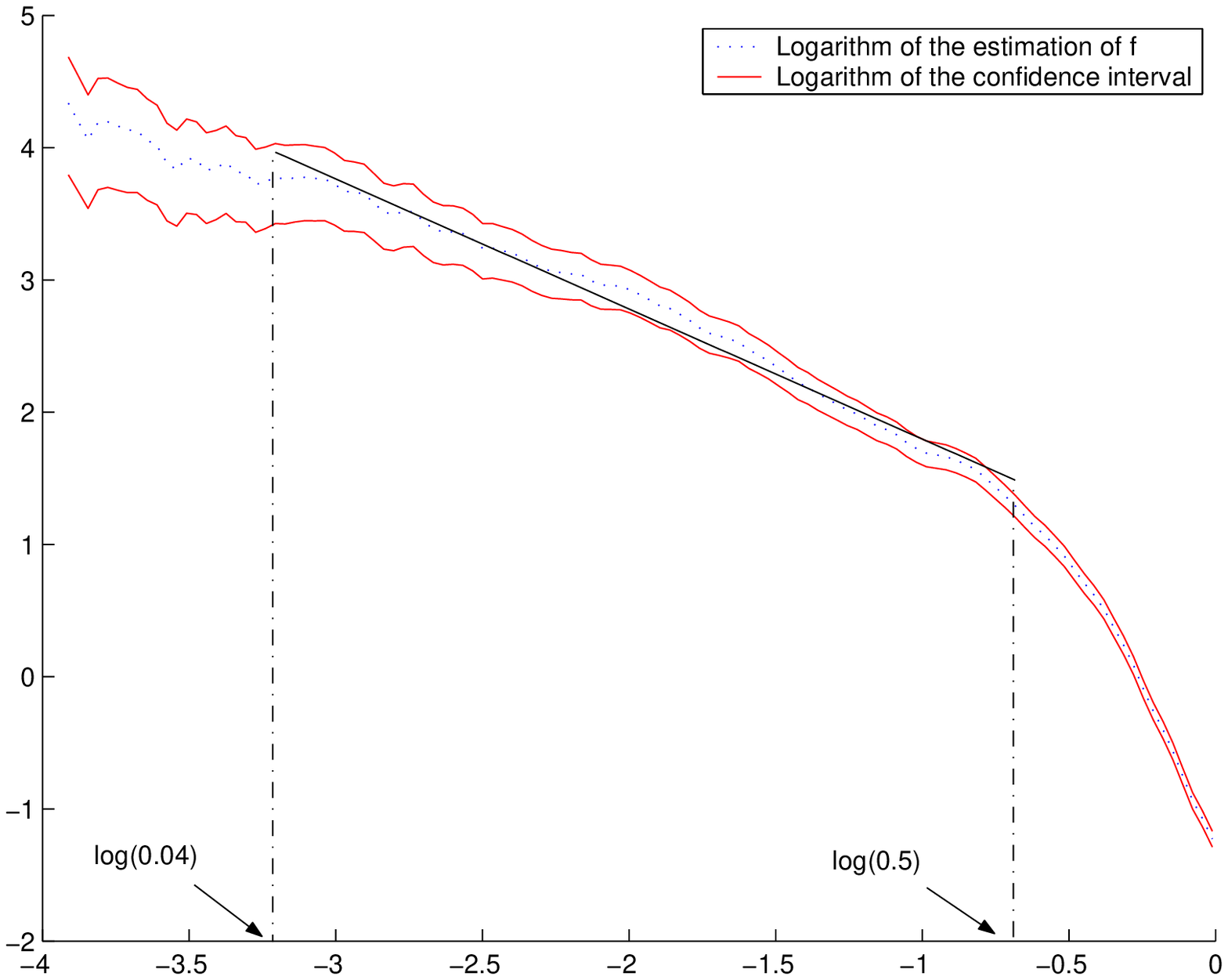}\]
\vspace{0.2cm}
\[
\epsfxsize 12cm \epsfysize 6.5cm \epsfbox{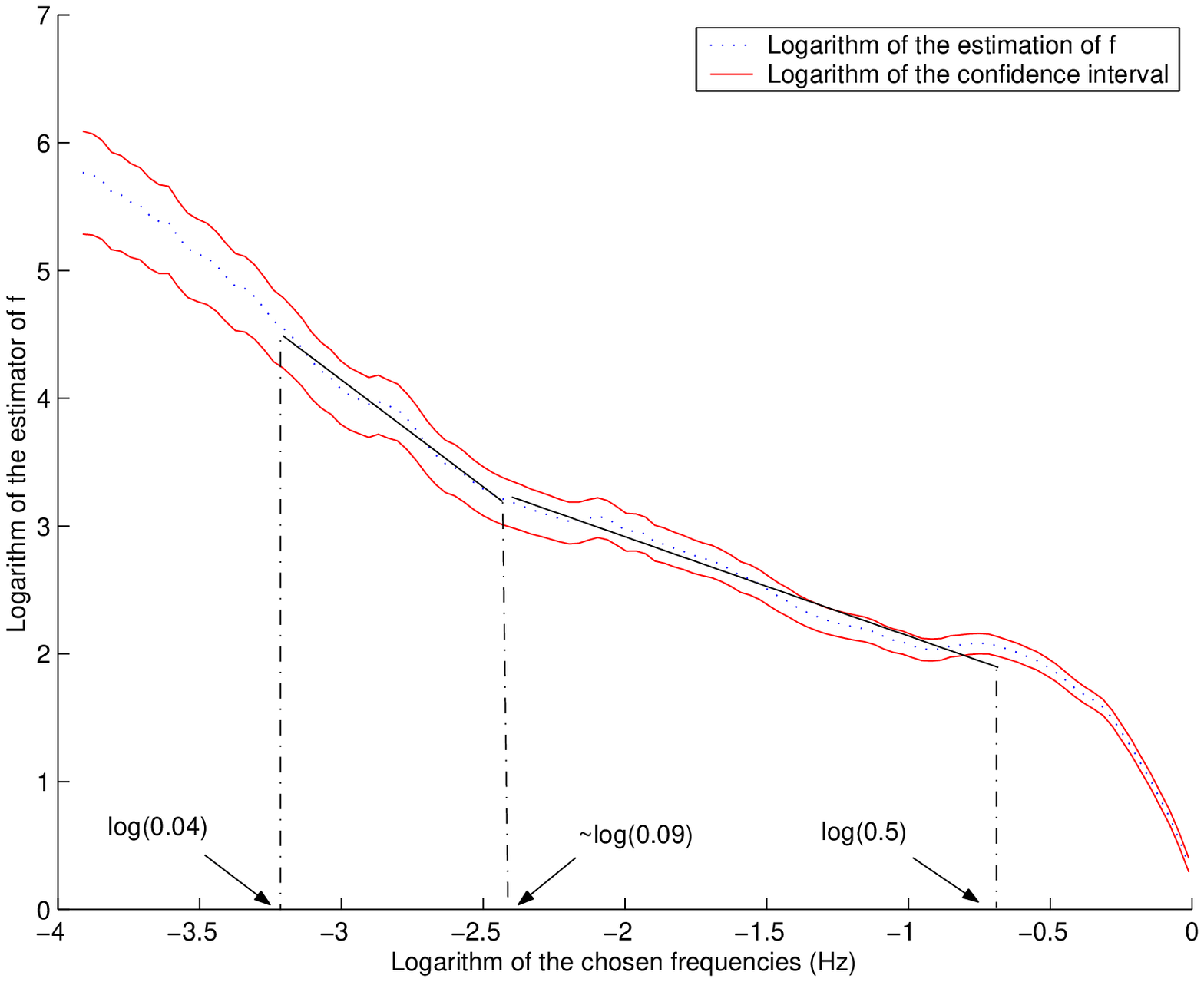}\]
\vspace{0.2cm}
\[
\epsfxsize 12cm \epsfysize 6.5cm \epsfbox{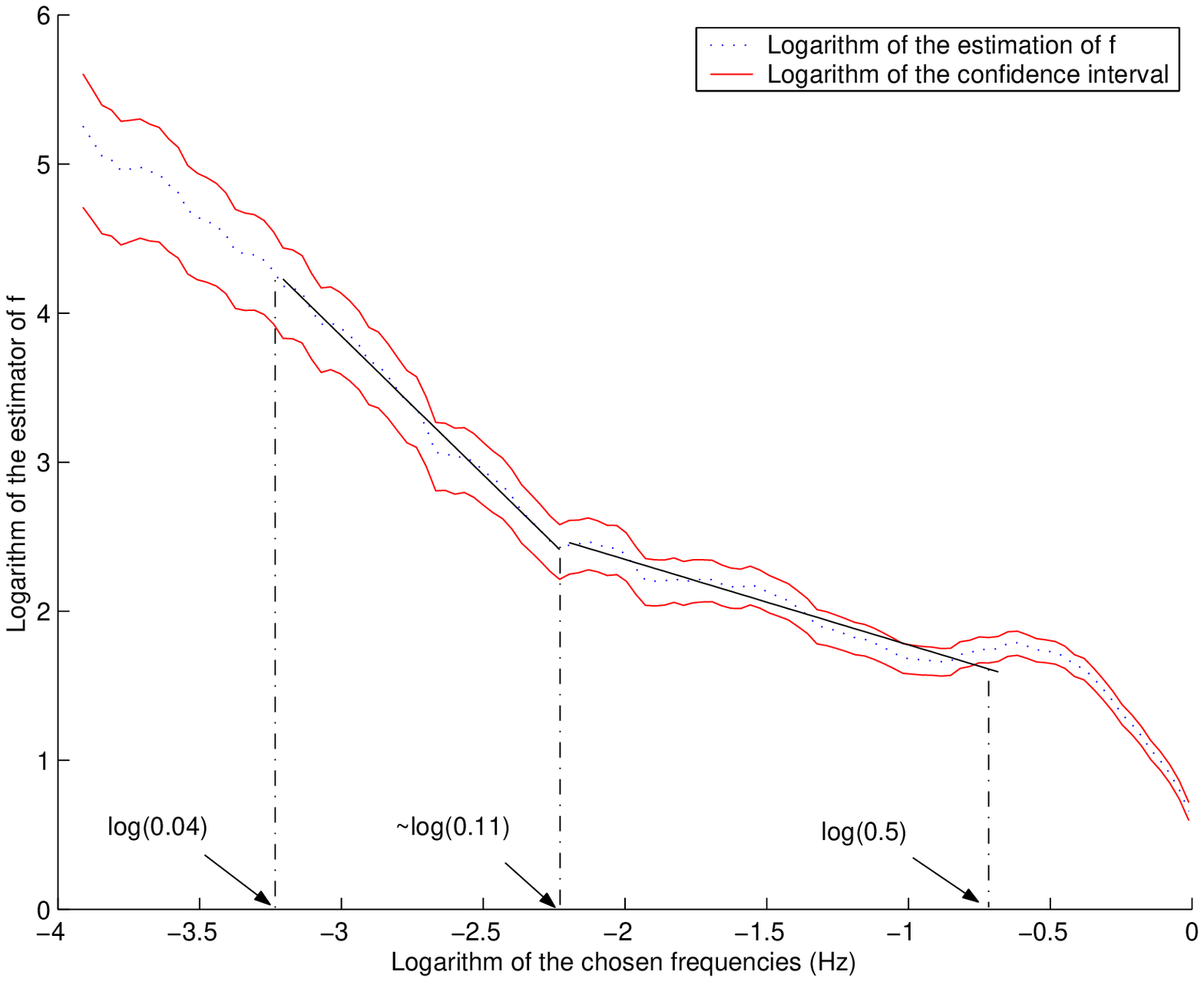}\]

\section*{Conclusion}
In this paper, we have constructed and studied an estimator of the spectral density $f$ of a Gaussian process from discrete sampling at instants $\big(t_i^{(n)}\big)_{i=0,\dots,n}$. One of the main novelties of this work is that the sampling scheme is random:
$\ds t_{i+1}^{(n)}-t_{i+1}^{(n)}= \delta_n L_i$ for a sequence of r.v. $L_i$ and a sampling step $\delta_n\to 0$. Under moment condition on the r.v. $L_i$, we have obtained a CLT for sample variance of wavelet coefficients. Then, by using the same wavelet-type technique with a bandwidth $\lambda_n\to \infty$, we have obtained a pointwise estimator of the spectral density at some frequency $\xi$. Under the conditions: $\lambda_n=o((n\delta_n)^{1/2})$ and $(n\delta_n)^{1/5}= o(\lambda_n)$, we have a CLT with a rate of convergence $(n\delta_n)^{-2/5}$.
The moment condition on r.v. $L_i$ is linked to the regularity index of the process $X$. This pointwise estimation of the spectral density is then applied to heartbeat time series. With this tool, we can observe variations of the spectral density according to the type of activity of human beings.

Three directions for future research have been opened by this work. Firstly, our results could be extended to non-Gaussian processes, {\it i.e.} when another measure replaces the Brownian measure in formula (\ref{repr:harmonizable}). Thus the expression of the spectral density estimator will be the same a non-Gaussian frame. However certain conditions will be necessary to establish CLTs similar to those obtained in theorem \ref{Theo:TCL:discret} and Proposition \ref{cvg:Ip}. Secondly, it could be interesting to make more investigations in order to compare our sampling scheme with Masry's conditions. Finally, from a practical point of view, we believe that our estimator has potential applications for other kinds of real data.
\par
\medskip
\noindent
{\bf Acknowledgments.} The authors are grateful to both the referees as well as editors for their very careful reading and many relevant
suggestions and corrections that strongly improve the content and
the form of the paper. We also thank A. Ayache, E. Masry and M. Wschebor for mathematical discussions. We thank A. Chamoux (Clermont-Ferrand Hospital) and the members of team UBIAE (INSERM and \'Evry Génopole) for discussions on heartbeat and physiological signal processing.

~\\
Corresponding author: Pierre, R.~Bertrand, 
{\it INRIA Saclay, 91893 Orsay Cedex, France, and Laboratoire de
Mathématiques, UMR CNRS 6620, Clermont Université, France, email:~
Pierre.Bertrand@math.univ-bpclermont.fr.} 
\section*{Appendix}
\subsection*{Proofs of useful lemmas and Proposition \ref{pte:d(a)}}
In the sequel, the following lemma will be useful:
\begin{lem}\label{maj:EXt:grand} Let $X$ be a Gaussian process defined by
(\ref{repr:harmonizable}) with a spectral density function $f$
satisfying (\ref{cond:f}) or by (\ref{stat}) with a spectral density satisfying (\ref{cond:f:stat}). Then there exists a constant $C_0>0$
such that  \ba \label{maj:cov1}
\left|\E\left(X(t_1) X(t_2)\right)\right| &\le& C_0 (1+ |t_1|) (1+
|t_2|)~~\mbox{for all $(t_1,t_2) \in \R^2$}. \ea
\end{lem}
\begin{dem} Firstly, let us consider $X$ defined by (\ref{repr:harmonizable}). For all $t\in \R$, by using (\ref{isometry:IS}), we have
\ban \E\big[ X^2(t)\big]&=& \int_{\R} \left|e^{it\xi} -
1\right|^2\, f(\xi)\, d\xi \,\le\, 2\int_0^1 |t \xi|^2\, f(\xi)\,
d\xi\,+\,8\int_1^{\infty} f(\xi)\, d\xi
\\
&\le& (2 t^2 + 8)\times \int_0^{\infty} \big(1\wedge |\xi|^2\big)\, f
(\xi) \, d \xi. \ean This implies $\displaystyle{\E\big(X(t)^2
\big) \le C_0\, (1+ |t|^2)}$ where $C_0 = 4\,\int_{\R}
\big(1\wedge  |\xi|^2\big)\, f (\xi) \, d \xi$. Then, by using
Cauchy-Schwartz inequality, one deduces (\ref{maj:cov1}).
Secondly, consider $X$ defined by (\ref{stat}), then (\ref{isometry:IS}) and (\ref{cond:f:stat}) imply that
\ban
\left|\E\left(X(t_1) X(t_2)\right)\right|
\,=\,
\left|\int_{\R} e^{i(t_1-t_2)\xi} f(\xi)\, d\xi\right|
&\le&
\int_{\R}  f(\xi)\, d\xi \,<\,\infty.
\ean
Therefore (\ref{maj:cov1}) is satisfied with $\ds C_0=\int_{\R}  f(\xi)\, d\xi$. This finishes the proof of Lemma \ref{maj:EXt:grand}.
\end{dem}
\begin{dem} [of Proposition \ref{pte:d(a)}]
We just give the proof for Gaussian processes defined by
(\ref{repr:harmonizable}) with a spectral density function $f$
satisfying (\ref{cond:f}). The key property, which explains that the same representation formula holds for
for Gaussian processes defined by  (\ref{repr:harmonizable}) or (\ref{stat}), is Formula (\ref{disparition:moins1}).
Moreover, since Condition (\ref{cond:f:stat}) implies Condition (\ref{cond:f}), all the convergence results remain valid under Condition (\ref{cond:f:stat}).

Let $X$ be defined by (\ref{repr:harmonizable}). Firstly, one can show that for all $a>0$ and $b\in \R$, $\E \left[\big|d_{X}(a,b)\big|^2\right]
< \infty$. This induces that $d_{X}(a,b)$ is well defined.
Indeed, since $X$ is a real valued process, one has
\ban \E \left[\big|d_{X}(a,b)\big|^2\right] &=& \frac 1 a
\int_{\R}\int_{\R} \psi\left(\frac{t_1-b}{a}\right)\,
\overline{\psi\left(\frac{t_2-b}{a}\right)} \E\left(X(t_1)
X(t_2)\right)dt_1dt_2
\\
&\le& \frac {C_0} a \int_{\R}\int_{\R}
\left|\psi\left(\frac{t_1-b}{a}\right)\right|
 \left|\psi\left(\frac{t_2-b}{a}\right)\right|\, (1+ |t_1|) (1+ |t_2|)dt_1dt_2
\\
& \leq &
a\, C_0 \left ( \int_{\R} \left|\psi(u)\right|\,(1+|b|+ |a\,u|)
du\right )^2< \infty,
\ean where we have used successively the bound (\ref{maj:cov1}),
the change of variable $u=(t - b)/a$  and the second condition of
Assumption W$(1, 1,0)$.
Next, one turns to the proof of
 the representation formula (\ref{repr:dX}). Firstly, recall that the stochastic of a complex valued function $g= g_1+i g_2$ against a complex Gaussian measure $W$ with real part $W_1$ and imaginary part $W_2$ is defined by
 $$\ds\int_{\R} g(x)\, dW(x) = \int_{\R} g_1(x)\, dW_1(x) -\int_{\R} g_2(x)\, dW_2(x)$$ and that $W_1$ and $W_2 $ are Wiener measures, see \cite[(7.2.8) p.326]{TaSa:1994}.
Now,  consider any
interval $[\alpha, A] \subset ]0, \infty[$, the
function $f$ is bounded on $[\alpha, A]$ 
 and $\displaystyle{\int_{\R}
\left|\psi\left(\frac{t-b}{a}\right)\right|\,dt = a\int_{\R}
\left|\psi(u)\right|\,du < \infty }$. Therefore, one can apply the
Fubini-type Theorem for stochastic integral (see \cite[Lemma 4.1,
p. 116]{IW}) to the two integrals 
corresponding to the real  and the imaginary part, then by summing up,   we get
\ban
\int_{\alpha}^{A} \left[\int_{\R}\left(e^{it\xi} - 1\right)\cdot\psi\left(\frac{t-b}{a}\right) \,dt\right]\,
f^{1/2}(\xi)\,dW(\xi) &=&
\int_{\R}\psi\left(\frac{t-b}{a}\right)\left[\int_{\alpha}^{A} \left(e^{it\xi} - 1\right)\cdot f^{1/2}(\xi)\,dW(\xi) \right]\,dt
\ean
Since $\ds\int_{\R}\psi(u)\,du =0$, for all couple $(a,b)\in
]0,\infty[ \times \R$, we have
\ba
\label{disparition:moins1}
\int_{\R}\left(e^{it\xi} - 1\right)\,\psi\left(\frac{t-b}{a}\right) \,dt
&=&
\int_{\R} e^{it\xi} \,\psi\left(\frac{t-b}{a}\right) \,dt
\,=\,
a\, e^{ib\xi}\, \overline{\widehat{\psi}}(a \xi).
\ea
Thus, for all couple $(a,b)\in ]0,\infty[ \times \R$ we have
\ba
\label{Fubini0}
a\, \int_{\alpha}^{A}e^{ib\xi}\, \overline{\widehat{\psi}}(a \xi)\,
f^{1/2}(\xi)\,dW(\xi) &=&
\int_{\R}\psi\left(\frac{t-b}{a}\right)\left[\int_{\alpha}^{A} \left(e^{it\xi} - 1\right)\cdot f^{1/2}(\xi)\,dW(\xi) \right]\,dt
\ea
From the one hand, the first condition of W$(1, 1, 0)$ and  (\ref{cond:f})  imply that
$\displaystyle{\,
\int_{0}^{\infty} \left|\widehat{\psi}(a
\xi)
\right|^2\, f(\xi)\,d\xi \,<\,
\infty.\,}$
 Therefore, one
can deduce that for any sequence of couples $(\alpha_n,\, A_n)$ converging to
$(0,\infty)$, the sequence $\displaystyle{\int_{\alpha_n}^{A_n}e^{iab\xi}\,
\overline{\hat{\psi}}(a
\xi)\, f^{1/2}(\xi)\,dW(\xi) }$ converges in $L^2(\Omega)$.  From
the other hand,
$$\displaystyle{\int_{\R}\psi\left(\frac{t-b}{a}\right)\left[\int_{\alpha_n}^{A_n}
\left(e^{it\xi} - 1\right)\cdot f^{1/2}(\xi)\,dW(\xi)
\right]\,dt}\,$$ converges also in $L^2(\Omega)$, because
$$\E \left|\int_{\R}\psi\left(\frac{t-b}{a}\right)\left[\int_{0}^{\infty}
\left(e^{it\xi} - 1\right)\cdot f^{1/2}(\xi)\,dW(\xi)
\right]\,dt\right|^2 < \infty.
$$
Indeed, firstly, by using Cauchy-Schwarz inequality, we get that
\ban
\left|\int_{0}^{\infty}
\left(e^{it_1\xi} - 1\right)\cdot \left(e^{-it_2\xi} - 1\right)\,f(\xi)\, d\xi\right|
&\le&
 \left(\int_{0}^{\infty}
\left|e^{it_1\xi} - 1\right|^2\,f(\xi)\,
d\xi\right)^{1/2}\times
 \left(\int_{0}^{\infty}
\left|e^{it_2\xi} - 1\right|^2\,f(\xi)\,
d\xi\right)^{1/2}
\\
&\hspace{-5.5cm}\le&\hspace{-3cm}
\left(\int_{0}^{1}
\left| t_1\xi \right|^2 f(\xi)\,
d\xi + 4\int_{1}^{\infty} f(\xi)\, d\xi\right)^{1/2}
\times
\left(\int_{0}^{1}
\left| t_2\xi \right|^2 f(\xi)\,
d\xi + 4\int_{1}^{\infty} f(\xi)\, d\xi\right)^{1/2}
\\
&\hspace{-5.5cm}\le&\hspace{-3cm}\left((4+ t_1^2)\,\int_{\R}
\left(1\wedge |\xi|^2\right)\cdot f (\xi)
\,
d \xi\right)^{1/2}\times \left((4+ t_2^2)\,\int_{\R}
\left(1\wedge |\xi|^2\right)\cdot f (\xi)
\,
d \xi\right)^{1/2}
\\
&\hspace{-5.5cm}\le&\hspace{-3cm} (2+t_1)\times (2+t_2)\times
\int_{\R}
\left(1\wedge |\xi|^2\right)\cdot f (\xi)
\,
d \xi.
\ean
Next, by using  the isometry property (\ref{isometry:IS}), we get the following upper bound
\ban
&&\hspace{-1cm}\E
\left|\int_{\R}\psi\left(\frac{t-b}{a}\right)\left[\int_{0}^{\infty}
\left(e^{it\xi} - 1\right)\cdot f^{1/2}(\xi)\,dW(\xi)
\right]\,dt\right|^2
\\
&=&
\int_{\R}\int_{\R}\psi\left(\frac{t_1-b}{a}\right)\,\overline{\psi\left(\frac{t_2-b}{a}\right)}
\left(\int_{0}^{\infty}
\left(e^{it_1\xi} - 1\right)\cdot \left(e^{-it_2\xi} - 1\right)\,f(\xi)\, d\xi\right)
\,dt_1\,dt_2
\\&\le&
\left(\int_{\R}
\left(1\wedge |\xi|^2\right)\cdot f (\xi)
\,
d\xi\right)
\times \int_{\R}\int_{\R}\psi\left(\frac{t_1-b}{a}\right)\,\overline{\psi\left(\frac{t_2-b}{a}\right)}
(2+t_1)\times (2+t_2)\,dt_1\,dt_2
\\
&=&
\left(\int_{\R}
\left(1\wedge |\xi|^2\right)\cdot f (\xi)
\,
d\xi\right)
\times
\left(\int_{\R}(2+t)\,\left|\psi\left(\frac{t-b}{a}\right)\right|\,dt\right)^2
\\
&<& \infty,
\ean
where the last bound follows from Condition (\ref{cond:f}) and
Condition W(1, 1, 0)). Eventually, one can pass to the limit in
(\ref{Fubini0}) which provides
\ban
a\, \int_{0}^{\infty}e^{ib\xi}\, \overline{\widehat{\psi}}(a \xi)\,
f^{1/2}(\xi)\,dW(\xi) &=&
\int_{\R}\psi\left(\frac{t-b}{a}\right)\left[\int_{0}^{\infty} \left(e^{it\xi} - 1\right)\cdot f^{1/2}(\xi)\,dW(\xi)
\right]\,dt.
\ean
But similar calculations would lead to the same result between the
bounds $-\infty$ and $0$. By adding the two integrals between
 $0$ and $\infty$ and between $-\infty$ and $0$, one can obtain
\ban
a\,\int_{\R} e^{ib\xi}\,  \overline{\hat{\psi}}(a \xi)\,
f^{1/2}(\xi)\,dW(\xi) &=&
\int_{\R}\psi\left(\frac{t-b}{a}\right)\left[\int_{\R} \left(e^{it\xi} - 1\right)\cdot f^{1/2}(\xi)\,dW(\xi)
\right]\,dt
\\
&=&
\int_{\R}\psi\left(\frac{t-b}{a}\right)\,X(t)\,dt
\ean
which implies (\ref{repr:dX}). Afterwards, formula (\ref{repr:dX})
implies that for all $a>0$, $b\in \R$, $d_{X}(a,b)$ is a Gaussian
centered random variable with variance $\mathcal{I}_{1}(a)$.
Moreover, for all $a>0$ and $(b_1,b_2) \in \R^2$, we have
\ban
\E\left(d_{X}(a,b_1)\,\overline{d_{X}(a,b_2)}\right)
&=& a\, \int_{\R} e^{ia(b_1-b_2)\xi}\,\left|\hat{\psi}(a
\xi)\right|^2\, f (\xi)\,d\xi
\ean
Thus for a given $a>0$,
$\E\left(d_{X}(a,b_1)\,d_{X}(a,b_2)\right)$ is only depending on
$(b_1-b_2)$ which induces that $(d_{X}(a,b))_{b\in \R}$ is a
stationary  process. This finishes the proof of Proposition
\ref{pte:d(a)}.
\end{dem}
\par\smallskip \noindent
From formula (\ref{repr:dX}), it is clear that  for all $(a_1,a_2)\in (0,\infty)^2$ and  for all
$(b_1,b_2)\in \R^2$,
\ba \nonumber
\E \big( d_{X}(a_1,b_1)\cdot\overline{d_{X}(a_2,b_2)}\big
)&=&\sqrt{a_1 a_2}\,\cdot \gamma(b_2-b_1,a_1,a_2)
\\
\nonumber\mathrm{where}\hspace{7cm}&&
\\
\label{def:gamma}
\gamma(\theta,a_1,a_2):&=& \int_{\R}
e^{i\theta\xi}\,\widehat{\psi}(a_1
\xi)\,\overline{\widehat{\psi}(a_2\xi)}\,f(\xi)\,d\xi.\hspace{4cm}
\ea
When $(a_1,a_2)$ are positive numbers, the function $\gamma$ and its first derivative with respect to $\theta$ can be bounded:
\begin{lem}\label{lem:gamma:a}
Under Assumption W$(1, 0, 1/2)$ and if $f$ satisfies (\ref{cond:f}) and Assumption F($H$) with $H>0$:
\begin{enumerate}
\item for all $(a_1,a_2)\in (0,\infty)^2$, there exists $C>0$ not depending on $\theta$ such
that,
$\displaystyle{\big|\gamma(\theta,a_1,a_2)\big |\,<C \, \big ( 1
\wedge |\theta|^{-1}\big ) }$ for all $\theta \in \R$.
\item the function $\gamma$ is derivable with respect to $\theta$
and for all $(a_1,a_2)\in (0,\infty)^2$, there exists $C>0$ not depending on $\theta$ such
that, $\ds |\gamma'(\theta,
 a_1,a_2)| :=\left |\frac{\partial\gamma}{\partial \theta}(\theta,
 a_1,a_2)\right| \le C' \, (1\wedge|\theta|^{-1})$ for all $\theta \in \R$.
\end{enumerate}
\end{lem}
\begin{dem}{\bf [of Lemma \ref{lem:gamma:a}]}
Firstly, from Assumption W$(1,0, 0)$ (induced by Assumption W$(1, 0, 1/2)$), there exists a constant $ c>0$
such that \ba \label{maj:hat:psi}
\big|\widehat{\psi}(\xi)\big|\le c\,\big (1 \wedge |\xi|^2\big
) \quad\mathrm{for\;all}~~ \xi \in \R. \ea Indeed, from one hand,
$\displaystyle{\big|\widehat{\psi}(\xi)\big|\le\|\psi\|_{L^1(\R)} <\infty }$. From
the other hand, $\psi  \in W(1,0, 0)$ implies that $\widehat{\psi}$
is twice continuously differentiable and $\widehat{\psi}(0) =
\widehat{\psi}'(0) = 0$. From Taylor-Lagrange Formula, for
all $\xi \in \R^*$, there exists $\xi_0\in \R$ with
$|\xi_0|\le |\xi|$ such that
$\displaystyle{\widehat{\psi}(\xi)\,=\, \frac 1 2\, \xi^2\times
\widehat{\psi}''(\xi_0)}$. This induces $\displaystyle
\big|\widehat{\psi}(\xi)\big| \le \frac 1 2\, |\xi|^2\times
\Big(\int_{\R} t^2\, |\psi(t)|\,dt\Big)$ providing the second
bound of (\ref{maj:hat:psi}).
\par
\bigskip
\noindent
Secondly, we  show the first item.  Inequality (\ref{maj:hat:psi}) implies
that \ban \int_{\R} \left|\overline{\widehat {\psi}}(a \xi) \right|^2\,
f(\xi)\,d\xi &\le& c^2 \,\Big ( \int_{|\xi|\leq 1} |a\xi|^4
f(\xi)\,d\xi +\int_{|\xi|> 1}  f(\xi)\,d\xi \Big ) \\
& \leq & c^2 \big (1 \vee a^4 \big )\int _{\R} \big (1 \wedge
\xi^2 \big ) f(\xi)\,d\xi <C, \ean with $C>0$ not depending
on $\theta$. From Cauchy-Schwarz Inequality,
$$
\gamma(\theta,a_1,a_2) \leq c^2 \big (1 \vee a_1^2 \big )\big (1
\vee a_2^2 \big )\, \int _{\R} \big (1 \wedge \xi^2 \big )
f(\xi)\,d\xi .
$$
Combined with (\ref{cond:f}), this means that $\gamma(\theta,a_1,a_2)$ is bounded by a constant.
Moreover, with $f(\omega_k^+)$ and $f(\omega_k^-)$ denoting
the right and left limits of $f$ at $\omega_k$, for
all $1\le k \le K-1$, $\theta\in \R^*$ and $(a_1,a_2)\in (0,\infty)^2$,
\ban
\int_{\omega_k}^{\omega_{k+1}} e^{i \theta \xi}
\overline{\widehat{\psi}}(a_1 \xi)\,\widehat{\psi}(a_2\xi)\,
f(\xi)\,d\xi&&
\\
&\hspace{-7.3cm}& \hspace{-3.5cm}= \frac 1 {i\theta}
\left(e^{i \theta\omega_{k+1}}f(\omega_{k+1}^-)\overline{\widehat{\psi}}(a_1\omega_{k+1})
\widehat{\psi}(a_2\omega_{k+1})-e^{i
\theta\omega_{k}}
f(\omega_{k}^+)\overline{\widehat{\psi}}(a_1\omega_{k})\widehat{\psi}(a_2\omega_{k})\right)
\\ &\hspace{-7.6cm}& \hspace{-3.5cm}-\int_{\omega_k}^{\omega_{k+1}}
\frac {e^{i \theta\xi}}
{i\theta}\left[f'(\xi)\overline{\widehat{\psi}}(a_1
\xi)\widehat{\psi}(a_2\xi)
+f(\xi)\Big(a_1\overline{\widehat{\psi'}}(a_1
\xi)\widehat{\psi}(a_2\xi)+a_2\overline{\widehat{\psi}}(a_1
\xi)\widehat{\psi'}(a_2\xi)\Big)\right]d\xi. \ean
The same
result remains in force for $k=0$ and $k=K$. Indeed, by using
(\ref{maj:hat:psi}) combined with Assumption F($H$), one deduces
that for all $\theta\in \R$ and $(a_1,a_2)\in (0,\infty)^2$,
\ban \lim_{\xi\to 0}e^{i
\theta\xi}\,f(\xi)\,\overline{\widehat{\psi}}(a_1\xi)\,
\widehat{\psi}(a_2\xi) =0 &\;\quad\mathrm{and}\quad\;&
\lim_{\xi\to \infty}e^{i
\theta\xi}\,f(\xi)\,\overline{\widehat{\psi}}(a_1\xi)\,
\widehat{\psi}(a_2\xi) =0. \ean
Moreover, since $f$ is an even function,
\ban
\int_{-\omega_{k+1}}^{-\omega_{k}} e^{i \theta \xi}
\overline{\widehat{\psi}}(a_1 \xi)\,\widehat{\psi}(a_2\xi)\,
f(\xi)\,dx&&
\\
&\hspace{-7.3cm}& \hspace{-3.5cm}= \frac 1 {i\theta}
\left(e^{-i \theta\omega_{k}}f(\omega_{k}^+)\overline{\widehat{\psi}}(-a_1\omega_{k})
\widehat{\psi}(-a_2\omega_{k})-e^{i
\theta\omega_{k+1}}
f(\omega_{k+1}^-)\overline{\widehat{\psi}}(a_1\omega_{k})\widehat{\psi}(a_2\omega_{k})\right)
\\ &\hspace{-7.6cm}& \hspace{-3.5cm}-\int_{-\omega_{k+1}}^{-\omega_{k}}
\frac {e^{i \theta\xi}}
{i\theta}\left[f'(\xi)\overline{\widehat{\psi}}(a_1
\xi)\widehat{\psi}(a_2\xi)
+f(\xi)\Big(a_1\overline{\widehat{\psi'}}(a_1
\xi)\widehat{\psi}(a_2\xi)+a_2\overline{\widehat{\psi}}(a_1
\xi)\widehat{\psi'}(a_2\xi)\Big)\right]d\xi. \ean
Thus, by summing up and using
Assumption F($H$), for all $ \theta\in \R$ and $(a_1,a_2)\in (0,\infty)^2$,
\ban
 \gamma(\theta,a_1,a_2)&=&  \sum_{k=0}^K\int_{\omega_k}^{\omega_{k+1}}  e^{i \theta \xi}
\overline{\widehat{\psi}}(a_1 \xi)\,\widehat{\psi}(a_2\xi)  \,
f(\xi)\,dx+ \sum_{k=0}^K\int_{-\omega_{k+1}}^{-\omega_{k}} e^{i \theta \xi}
\overline{\widehat{\psi}}(a_1 \xi)\,\widehat{\psi}(a_2\xi)\,
f(\xi)\,dx\\
&=& -\frac {1}
{i\theta} \sum_{k=1}^{K} \Big ( e^{i \theta\omega_{k}}
\overline{\widehat{\psi}}(a_1
\omega_{k})\widehat{\psi}(a_2\omega_{k})-e^{-i \theta\omega_{k}}
\overline{\widehat{\psi}}(-a_1
\omega_{k})\widehat{\psi}(-a_2\omega_{k})\Big )\big(
f(\omega_{k}^+)-f(\omega_{k}^-)\big)\\
&& + \frac {1} {i\theta} \int_{\R}
 e^{i \theta\xi} \Big [f'(\xi)\overline{\widehat{\psi}}(a_1
\xi)\,\widehat{\psi}(a_2\xi)
+f(\xi)\Big(a_1\overline{\widehat{\psi'}}(a_1
\xi)\widehat{\psi}(a_2\xi)+a_2\overline{\widehat{\psi}}(a_1
\xi)\widehat{\psi'}(a_2\xi)\Big)\Big ] d\xi\qquad \ean since
the integral of the r.h.s. of the previous equality is well
defined.  Then,
\ban
\left|\gamma(\theta,a_1,a_2)\right| &\le& \frac {1} {|\theta|}
\Big ( 2c \, \sum_{k=1}^{K}\big  | f(\omega_{k}^+)
\,-\,f(\omega_{k}^-)\big|\\
&+& \int_{\R}
\Big [\big|f'(\xi)\overline{\widehat{\psi}}(a_1
\xi)\widehat{\psi}(a_2\xi)\big|
+|f(\xi)|\big(|a_1|\big|\overline{\widehat{\psi'}}(a_1
\xi)\widehat{\psi}(a_2\xi)\big|+|a_2|\,\big|\overline{\widehat{\psi}}(a_1
\xi)\widehat{\psi'}(a_2\xi)\big|\big)\Big]d\xi \Big). \ean It
remains to show the convergence of the previous integral. Using the same trick
as in Formula
(\ref{maj:hat:psi}), under Assumption W$(1,0, 0)$, $\big|\widehat{\psi'}(\xi)\big|\le
c'\,
\big (1\wedge |\xi| \big )$ with $c'$ not depending on $\xi$. So, for all $(a_1,a_2)\in (0,\infty)^2$
\ban
&&\hspace{-0.6cm}\int_{\R}
\Big [\big|f'(\xi)\overline{\widehat{\psi}}(a_1
\xi)\widehat{\psi}(a_2\xi)\big|
+|f(\xi)|\big(|a_1|\big|\overline{\widehat{\psi'}}(a_1
\xi)\widehat{\psi}(a_2\xi)\big|+|a_2|\,\big|\overline{\widehat{\psi}}(a_1
\xi)\widehat{\psi'}(a_2\xi)\big|\big)\Big]d\xi \\
&& \hspace{+0.2cm} \le C_1 \, a_1^2a_2^2\,\int_{|\xi|\le 1}
|f'(\xi)|\xi^4 + 2|f(\xi)\xi^3|\,d\xi +C_2\,
\int_{|\xi|>1} |f'(\xi)| +(|a_1|+|a_2|)|f(\xi)|\,d\xi \\
&&\hspace{0.2cm} \leq C\, \int_{\R}\Big [\big (1\wedge |
\xi|^4\big)\cdot|f'(\xi)| + \big(1\wedge |
\xi|^3\big)\cdot |f(\xi)|\Big ] \,d\xi \ean
where $C>0$ depends on $c$, $c'$, $a_1$ and $a_2$. But since $\big (1\wedge |
\xi|^4\big) \leq \big (1\wedge |
\xi|^3\big)$ and $\big (1\wedge |
\xi|^3\big) \leq \big (1\wedge |
\xi|^2\big)$, then from Assumption F($H$), $\int_{\R} \big (1\wedge |
\xi|^4\big)\cdot|f'(\xi)|<\infty$ and $\int_{\R} \big (1\wedge |
\xi|^3\big)\cdot|f(\xi)|<\infty$ and this completes the proof of the first item.
Eventually, one proves the second item. The differentiability is
obvious and  \ban \gamma'(\theta,a_1,a_2) &=& i \int_{\R}
e^{i\theta\xi}\,\xi\,\overline{\widehat{\psi}}(a_1\xi)\,\widehat{\psi}(a_2\xi)\,f(\xi)\,d\xi.
\ean Assumption W$(1, 0, 1/2)$ implies that for all $a>0$,~~
$\displaystyle{|a\xi|^{1/2} |\widehat{\psi}(a \xi)|\le\,
C_{\psi}}$ for all $\xi \in \R$. Combined with (\ref{maj:hat:psi}), this induces that for all
$a>0$ and $\theta
\in\R$,
\ban
\left|\varphi'(\theta,a_1,a_2)\right|
&\le&\int_{\R}|\xi|\,|\widehat{\psi}(a_1
\xi)|\,|\widehat{\psi}(a_2\xi)|\,f(\xi)\,d\xi
\\
&\le& c^2 (a_1\, a_2)^2 \int_{|\xi|\le 1}
|\xi|^5\,f(\xi)\,d\xi\,+\,\frac{C^2_{\psi}}{\sqrt{a_1\, a_2}}
\int_{|\xi|> 1}   f(\xi)\,d\xi
\\
&\le& C, \ean with $C>0$ not depending on $\theta$. Using the same arguments as above, for all
$ \theta\in \R^*$ and $(a_1,a_2)\in (0,\infty)^2$,
\ban
\gamma'(\theta,a_1,a_2)
&=& -\frac {1} {\theta}
\sum_{k=1}^{K} \Big (e^{i \theta\omega_{k}}\omega_{k}
\overline{\widehat{\psi}}(a_1
\omega_{k})\widehat{\psi}(a_2\omega_{k})+ \\
&&\hspace{1cm} +e^{-i
\theta\omega_{k}}\omega_{k} \overline{\widehat{\psi}}(-a_1
\omega_{k})\widehat{\psi}(-a_2\omega_{k})\Big )\big(
f(\omega_{k}^+)-f(\omega_{k}^-)\big) \\
&-&\frac {1} {\theta} \int_{\R}
 e^{i \theta\xi} \Big [f(\xi)
\overline{\widehat{\psi}}(a_1\xi)\widehat{\psi}(a_2\xi)+ \xi
f'(\xi)\overline{\widehat{\psi}}(a_1 \xi)\widehat{\psi}(a_2\xi)
+\\
&& \hspace{1cm}+
\xi f(\xi)\Big ((a_1\,\overline{\widehat{\psi'}}(a_1
\xi)\widehat{\psi}(a_2\xi)+a_2\overline{\widehat{\psi}}(a_1
\xi)\widehat{\psi'}(a_2\xi)\Big )\Big ]d\xi
\ean
and therefore $\ds \left|\gamma'(\theta,a_1,a_2)\right|  \le
 \frac C {|\theta|},$
with $C>0$ not depending on
$\theta$. This finishes the proof.
\end{dem}
\vspace{-1cm}
\subsection*{Asymptotic behavior of sample variances of wavelet coefficients for continuous time processes}
Since $ \mathcal{I}_{1}(a)$ is
obviously defined from $\big | d_{X}(a ,b)\big |^2$, we begin with the study of
\ba \label{def1:IN(a)} I_n(a) &:=& \frac 1 {n+1} \,  \sum_{k=0}^{n} \big |
d_{X}(a ,c_k)\big |^2,\quad \mbox{for $a>0$ and $n\in \N^*$}. \ea
For $ n\in
\N^*$ and $a>0$, define also:
\ba
\label{S_n} \hspace{1cm} S_n^{2}(a): =\frac {2\,a^2}{(n+1)^2}\,\sum_{k=0}^{n}\sum_{\ell=0}^{n} \left| \int_{\R} e^{i\big
(\E(c_k-c_{\ell})\big )\xi}\,\big |
\widehat{\psi}(a \xi)\big |^2 f(\xi)\,d\xi \right|^2.
\ea
Note that $S_n$ depends on $\E (c_k-c_{\ell})$ and therefore its formula is valid when $(c_k)$ are r.v. However, we begin by proving the following proposition in the case of deterministic $(c_k)$ and of course $\E (c_k-c_{\ell})$ can be replaced by $c_k-c_{\ell}$ in (\ref{S_n}).
\begin{Prop}\label{TCL:ponctuel}
Let $X$ be a Gaussian process defined by (\ref{repr:harmonizable}) or (\ref{stat})
with a spectral density $f$ satisfying (\ref{cond:f}),
$\psi$ satisfying Assumption W$(1, 1, 1/2)$. Then if
$(c_k)_k$ is a family of real numbers such that $c_{1}< c_2<\ldots<
c_{n}$, $~~  n\,
\max_{1 \leq k \leq n}\{c_{k+1} -c_k\} \limiten \infty~$ and there exists a constant $C''>0$ such that for all  $n\in \N^*$
$$
\max_{1 \leq k \leq n}\{c_{k+1}
-c_k\} \leq C''\min_{1 \leq k \leq n}\{c_{k+1}-c_k\}<\infty
$$
then
$\forall a>0$, \ba \label{TCL1} \frac 1 {S_n(a)}\,  \big (I_n(a)
- \mathcal{I}_{1}(a)\big ) \limiteloin \mathcal{N}(0,1). \ea
Moreover, there exist $0<C_m<C_M$ not depending on
$n$ such that  $\forall n \in N^*$,
\begin{eqnarray}\label{inegSn}
C_m \leq S_n(a)\, \big(n \, \max_{1 \leq k \leq n}\{c_{k+1} -c_k\} \big)^{1/2} \leq C_M.
\end{eqnarray}
\end{Prop}
The proof of Proposition \ref{TCL:ponctuel} relies on Lemma \ref{lem:gamma:a} and the
following Lemma which is a Lindeberg CLT
(see a proof in Istas and Lang, 1997):
\begin{lem}\label{lem1}
Let $(Y_{N,i})_{1\leq i \leq N,\, N\in \N^*}$ be a triangular array of
zero-mean Gaussian r.v. Define $S_N^2 := \var (V_N)$ with $V_N := \sum_{
i=1}^{N} Y^2_{N,i}$ and $\displaystyle \beta_N :=\max_{1\le i\le
N}\sum_{j=1}^{N}|\cov\left(Y_{N,i}, Y_{N,j}\right)|$.
If $\ds \lim_{N\to \infty} \frac {\beta_N}{S_N} = 0$, then $S_N^{-1}
(V_N-\E(V_N))$ converges weakly to a standard Gaussian random
variable.
\end{lem}
\begin{dem}[{\bf of Proposition \ref{TCL:ponctuel}}]\\
Consider $\displaystyle{Y_{n,i} =
(n+1)^{-1/2}\, d_X(a,c_{i})}$ for
$i=0,\ldots,n$ and \ban
\left \{ \begin{array}{l}  \beta_n
= (n+1)^{-1} \, \max_{1 \leq i \leq  n} \Big \{ \sum_{j=0}^{n} |\cov\big( d_X(a,
c_{i}), d_X(a,c_{j})\big )|\Big \}, \\
S_n^2 =(n+1)^{-2} \, \sum_{i=0}^{n}\sum_{j=0}^{n}
\cov\big( d^2_X(a,c_{i}), d^2_X(a,c_{j})\big)
\end{array}\right . . \ean
But, by using Formula (\ref{def:gamma}),
$\forall (a,a_1,a_2) \in (0,\infty[^3$, $(b_1,b_2)\in \R^2$ \ban
\cov\big( d_X(a,b_1), d_X(a, b_2) \big) &=& a\,
\gamma\big(b_1-b_2, a, a\big)\\
\cov\big(d_X^2(a_1,b_1),\, d_X^2(a_2,b_2)\big )
&=& 2\,\left(a_1 a_2\right)\gamma^2(b_1-b_2,
a_1,a_2), \ean since variables $d_X(a,b)$ are zero-mean Gaussian
r.v. Therefore,
\ban \left \{ \begin{array}{l} \beta_n
= a\,(n+1)^{-1} \,\max_{0 \leq i\leq  n}\Big \{ \sum_{j=0}^{n}
\big |\gamma\big(c_{i}-c_{j}, a, a\big)\big | \Big \} \\
S_n^2 =2\,a^2\,(n+1)^{-2} \,\sum_{i=0}^{n}\sum_{j=0}^{n}
\gamma^2\big(c_{i}-c_{j}, a, a\big)\end{array}\right .  . \ean
Let $p$ and $q$ be such that $
1/p + 1/q = 1$ with $(p,q) \in (1,\infty)^2$. Then the H\"older
Inequality implies that \ban \beta_n &\le& C\,a\,  n^{1/q-1}\times
\max_{0 \leq i\leq  n}\Big\{\Big ( \sum_{j=0}^{n}\,
\Big|\gamma\big(c_{i}-c_{j}, a, a\big)\Big|^p\Big
)^{1/p} \Big \}. \ean
Lemma \ref{lem:gamma:a} 1) implies that for
every $1 \leq i\leq  n$, for $n$ large enough, \ba \nonumber
\sum_{j=0}^{n}\, \Big|\gamma\big(c_{i}-c_{j}, a,
a\big)\Big|^p & \le& C \, \Big (
\#\left\{0 \leq j \leq  n,\, \big|c_{i}
-c_{j}\big| \le 1\right\}+\sum_{j=0}^{n}
\big|c_{i} -c_{j}\big|^{-p}\1_{\big|c_{i} -c_{j}\big|>1} \Big )\\
\nonumber && \hspace{-3.5cm} \leq C \Big (
\#\big\{0 \leq j \leq  n,\, |i-j|\min_{0 \leq k\leq  n-1}\big |c_{k+1} -c_{k}\big | \le 1\big\}\\
&& \hspace{-0.9cm}+\sum_{j=0}^{n}
 \Big [
|i-j|\min_{0 \leq k\leq  n-1}\big |c_{k+1} -c_{k}\big |\Big
]^{-p }\1_{|i-j|\max_{0 \leq k\leq  n-1}\big |c_{k+1} -c_{k}\big |>1} \Big )\\
\label{borne_betan}  && \hspace{-3.5cm} \leq 2 \, C \,
\big(\min_{0 \leq k\leq  n-1}\big |c_{k+1} -c_{k}\big
|\big)^{-1} \Big(1 +\big(\min_{0 \leq k\leq  n-1}\big |c_{k+1}
-c_{k}\big |\big )^{1-p} \hspace{-1.7cm} \sum _{\ell \geq
(\max_{0 \leq k\leq  n-1} |c_{k+1} -c_{k}
|)^{-1}}\hspace{-1.2cm} |\ell|^{-p} \Big). \ea
Since $p>1$,
$\displaystyle \sum _{\ell=1}^{\infty} |\ell|^{-p}<\infty $ is
finite and thus $$ \sum _{\ell \geq (\max_{0 \leq k\leq  n-1}\big
|c_{k+1} -c_{k}\big |)^{-1}}\hspace{-1.2cm}
|\ell|^{-p} \leq \frac 1 {p-1} \, \big (\max_{0 \leq k\leq  n-1}\big
|c_{k+1} -c_{k}\big |\big)^{p-1}\leq \frac {C''} {p-1} \, \big (\min_{0 \leq k\leq  n-1}\big
|c_{k+1} -c_{k}\big |\big)^{p-1},$$
since by definition $\max_{1 \leq k \leq n}\{c_{k+1}
-c_k\} \leq C''\min_{1 \leq k \leq n}\{c_{k+1}-c_k\}<\infty$. Therefore,
\ba
\label{bound_betan} \beta_n
 &\le&
C\,a\,\left\{n\times\min_{0 \leq k\leq  n-1}\big |c_{k+1}
-c_{k}\big |\right\}^{-1/p}, \ea with $C>0$ not depending on $n$. 
Now, a lower bound for
$S_n^2$ is required. For all $a>0$,  $\theta \in \R \mapsto\gamma(\theta, a,
a)$ is a continuous map and  $\gamma(0, a, a) =
\int_{\R} \big|\widehat{\psi}(a\xi)\big|^2\,f(\xi)\,d\xi \,>\,0
$.
Therefore, for all $a >0$,
there exists $\theta_a>0$ such that $\ds\gamma(\theta, a, a)\ge
\frac 1 2\, \gamma(0, a, a)$ when $|\theta|\le
\theta_a$. Then,
\ba \nonumber S_n^2(a)\hspace{-2mm} &\ge& \hspace{-2mm}C'_1\,a^2\,n^{-2} \,
\gamma^2(0, a,a)\, \#\big \{0 \leq i,j \leq n,\,
\big|c_{i}
-c_{j}\big|\le\theta_a\big \} \\
\nonumber \hspace{-2mm}&\ge&\hspace{-2mm} C'_1\,a^2\,n^{-2}\,\gamma^2(0, a,a)\,
\# \Big \{0 \leq i,j \leq n,\,|i-j|\max_{0 \leq k \leq n-1} \big|c_{k+1}
-c_{k}\big|\le\theta_a\Big\} \\
\nonumber &\geq  & \hspace{-2mm}  C'_1\, a^2\,n^{-2} \gamma^2(0, a,a)
\,\frac  n 2 \, \Big (  \big (\theta_a\,\max_{0 \leq k \leq n-1}\big|c_{k+1}
-c_{k}\big|\big)^{-1} \wedge (n-1) \Big ), \ea
since for $\mu>0$, $\displaystyle \# \Big \{0 \leq i,j \leq n,\,|i-j|\le \mu \Big\}=2\sum_{k=0}^{[\mu] \wedge (n-1)} n-k  \geq 2 \big ( [\mu] \wedge (n-1) \big ) \,\frac n 2 $. Thus, since from assumptions $\displaystyle n\,\max_{0 \leq k \leq n-1} \big|c_{k+1}
-c_{k}\big| \limiten \infty$, there exists $C_M>0$ such that for $n$ large enough:
\ba
\label{bound_Sn}S_n^2(a) \geq   C_M \, \big (n \, \max_{0 \leq k \leq n-1}\big|c_{k+1}
-c_{k}\big|\big)^{-1}.
\ea
Now, from (\ref{bound_betan}) and  (\ref{bound_Sn}),
 \ban
\frac{\beta_n}{S_n}\le C \,  n^{1/2-1/p} \big( \max_{0 \leq k \leq n-1}
\big|c_{k+1} -c_{k}\big| \big)^{1/2} \, \big (
\min_{0 \leq k \leq n-1}\big |c_{k+1} -c_{k}\big |\big )
^{-1/p}. \ean
Therefore $\displaystyle{ \beta_n/S_n \le C\,
\big(n\,\max_{0 \leq k \leq n-1} \big|c_{k+1}
-c_{k}\big|\big)^{1/2-1/p}}$ with
$C>0$.
Next for any $p\in (1,2)$,\\
$\displaystyle{\lim_{n\to \infty}\big(n\,\max_{0 \leq k \leq n-1}
\big|c_{k+1} -c_{k}\big|\big)^{1/2-1/p} =0}$
and thus, $\displaystyle{\lim_{n\to
\infty}
\beta_n/S_n =0}$ and assumptions of Lemma
\ref{lem1} are fulfilled.\\
~\\
Finally, from (\ref{bound_Sn}),
$\ds S^2_n(a) \ge C_M \, \big(n  \max_{0 \leq k \leq n-1}\big
|c_{k+1} -c_{k}\big |\big)^{-1}$ with $C_M>0$ for $n$ large enough. Moreover, using
the bound (\ref{borne_betan}) for $p=2$ and the lines after (\ref{borne_betan}),
\ban
\sum_{j=0}^{n} \,
\gamma^2\big(c_{i}-c_{j}, a, a\big)  & \leq &  C \,
\big(\min_{0 \leq k \leq n-1}\big |c_{k+1} -c_{k}\big
|\big)^{-1}\\
\Longrightarrow ~S_n^2(a) \leq C'\,a^2\,n^{-2}\sum_{i=0}^{n}
\sum_{j=0}^{n} \gamma^2\big(c_{i}-c_{j}, a,
a\big)
 & \leq & C_m\,   \big (n \, \max_{0 \leq k \leq n-1}\big
|c_{k+1} -c_{k}\big |\big )^{-1}. \ean  Therefore, inequalities (\ref{inegSn}) are proved.
\end{dem}
\begin{Prop}\label{TCL:ponctuel2}
Let $X$ be a Gaussian process defined by (\ref{repr:harmonizable}) or (\ref{stat})
with a spectral density $f$ satisfying (\ref{cond:f}),
$\psi$ satisfying Assumption W$(1,1,r)$ with $r>1/2$. Then if
$(c_k)_k$ is a family of r.v. independent to ${\cal F}_X$ such that $c_{k}=c_0+\frac k n (c_n-c_0)$, with
$$ n^{\frac 1 {2r}-1} \, \E (c_n-c_0) \limiten 0,\quad  \frac {\E (c_n-c_0)}{\log n} \limiten \infty ~~\mbox{and}~~  \var (c_n-c_0)\limiten 0.
$$
Then (\ref{TCL1}) holds with
\ba \label{SN}
\lim_{n \to \infty} \big (\E (c_n-c_0) \big )\, S^2_n(a)= 4 \pi \,
a^2   \, \int_{\R}\big |
\widehat{\psi}(a z)\big |^4 f^2(z)\, d z.
\ea
\end{Prop}
\begin{rem}\label{casck}
For $(c_k)_k$ satisfying (\ref{ck}), under Assumption S($2$), Proposition \ref{TCL:ponctuel2} holds when
$n^{1/2} \, \delta_n \limiten 0$ because $\E |T_n - \E T_n|^{2} \leq n\,  \delta_n^2 \, \max_{1\leq k \leq n} \E L_k$.
\end{rem}
\begin{dem}[{\bf of Proposition \ref{TCL:ponctuel2}}]
The sequence of r.v. $(c_{k})_{0\leq k \leq n}$ is  independent to ${\cal F}_X$. Therefore, $(d_X(a,c_k))_{0\leq k \leq n}$ as the same distribution than $(d_X(a,c_k-c_0))_{0\leq k \leq n}$. Indeed for a sequence of deterministic real numbers $(c_k)_{0\leq k \leq n}$, $(d_X(a,c_k))_{0\leq k \leq n}$ is a stationary sequence and after $(d_X(a,c_k))_{0\leq k \leq n}$ has the same distribution than $(d_X(a,c_k-b_0))_{0\leq k \leq n}$. Next, conditionally to the $\sigma$-algebra $\sigma((c_k)_{0\leq k \leq n})$, $(d_X(a,c_k))_{0\leq k \leq n}$ as the same distribution than $(d_X(a,c_k-c_0))_{0\leq k \leq n}$. Finally, since $\sigma((c_k)_{0\leq k \leq n})$ and ${\cal F}_X$ are independent, $(d_X(a,c_k))_{0\leq k \leq n}$ as the same distribution than $(d_X(a,c_k-c_0))_{0\leq k \leq n}$. \\
Now, we can only consider the case: $c_k=k \tau_n/n$ with $\tau_n:=c_n-c_0$.
Define \ban I'_n(a) &:=& \frac 1 {n} \, \sum_{k=0}^{n} d_{X}^2\big (a ,\E c_k\big).
\ean It is clear that $(\E c_k)_{1 \leq k \leq n}$
is a deterministic sequence.
Thus  \ba
\label{TLC_stoch}
\frac{I'_n(a) - \mathcal{I}_{1}(a) }{S_n(a)}\limiteloin
\mathcal{N}(0,1). \ea
Nowadays, one has to check that the error $I'_n(a) -I_n(a)$ is negligible with respect to $S_n(a)$ in norm $L^2(\Omega)$.
But
$$
S_n(a)\ge C_M \cdot\big(n\max_{0\leq k\leq  n-1}\big|\E \big (c_{k+1}
-c_{k}\big )\big|\big)^{-1/2} \ge C \times  \big ( \E\tau_n\big)^{-1/2}.$$
Therefore, it suffices to prove that
$\ds
\lim_{n\to \infty} \E\tau_n \times \E \big [ (I'_n(a)-I_n(a))^2\big ] = 0.\,
$
Since the r.v.
$c_k$ are independent on ${\cal F}_X$, one gets
\ban
\E \big [ (I'_n(a)-I_n(a))^2\big ]   & =&   \E \Big [\E
\big [ (I'_n(a)-I_n(a))^2~\big | ~{\cal F}_X\big ] \Big ] \\
&& \hspace{-2.1cm} =  \frac 1 {(n+1)^{2}}   \sum_{k,k'=0}^{n}   \E \Big [
\E \big [ \big ( d_{X}^2\big (a ,\E c_k\big)-d_{X}^2\big
(a ,c_k\big)\big )  \big ( d_{X}^2\big (a
,\E c_{k'}\big)-d_{X}^2\big (a ,c_{k'}\big)\big )
\big | {\cal F}_X \big ]\Big ]\\
&&\hspace{-2.1cm} = \frac {2 a^2} {(n+1)^{2}}   \sum_{k,k'=0}^{n}  \E \Big [
\gamma^2 \big (\E c_k-\E c_{k'},a,a\big )-
\gamma^2\big (\E c_k-c_{k'},a,a\big )
\\
&& \hspace{1.3cm} -\gamma^2\big
(c_k-\E c_{k'},a,a\big ) +\gamma^2\big
(c_k-c_{k'},a,a\big ) \Big ]. \ean
Next, from Taylor
expansions,
\ban \gamma^2\big
(\E c_k  - c_{k'},a,a\big ) &=& \gamma^2
\big(\E c_k - \E  c_{k'},a,a\big ) +  2 \big ( \E c_{k'} - c_{k'} \big )\times\cdots\\
&&\hspace{-1.8cm}\int_0^1 \gamma\big
(\E c_k - \E  c_{k'}+ \lambda\,(\E  c_{k'}- c_{k'}),a,a \big )\gamma'\big
(\E c_k - \E  c_{k'}+ \lambda\,(\E  c_{k'}- c_{k'}),a,a \big )
 d\lambda\\
\gamma^2\big (c_k  - \E c_{k'},a,a\big ) &=&  \gamma^2
\big(c_k - c_{k'},a,a\big ) +  2  \big (
c_{k'} - \E c_{k'} \big )\times\cdots
\\
&&\hspace{-1.8cm}\int_0^1 \gamma\big
(c_k -  c_{k'}+ \lambda\,(\E  c_{k'}- c_{k'}),a,a \big )\,\gamma'\big
( c_k -  c_{k'}+ \lambda\,(\E  c_{k'}- c_{k'}),a,a \big )\,
 d\lambda.\ean
From Lemma \ref{lem:gamma:a},
 $\exists C>0$ such that $\ds \big|\gamma(\theta,a,a
)\gamma'(\theta,a,a )\big|\le   C\times \big(1\wedge\theta^{-2}\big)$ for all $\theta\in \R$. One
can deduce that
\ban
\Big|
\gamma^2 \big (\E c_k-\E c_{k'},a,a\big )-
\gamma^2\big (\E c_k-c_{k'},a,a\big )
- \gamma^2\big
(c_k-\E c_{k'},a,a\big ) +\gamma^2\big
(c_k-c_{k'},a,a\big )\Big| &&
\\
&&\hspace{-11cm}\le \, C \,\big | c_{k'}-\E c_{k'}\big|
\times \int_0^1 \big ( 1\wedge\theta _{1,k,k'}^{-2}(\lambda)\big) +
\big ( 1\wedge\theta _{2,k,k'}^{-2}(\lambda)\big)\,d\lambda
\ean
with
$\theta_{1,k,k'}(\lambda) = \E (c_k -  c_{k'})+ \lambda(\E  c_{k'}- c_{k'})$ and
$\theta_{2,k,k'}(\lambda) = c_k -  c_{k'}+ \lambda(\E  c_{k'}- c_{k'})$.
~\\
Then,
\begin{eqnarray} \label{erreur}
\E \big [ (I'_n(a)-I_n(a))^2\big ]\leq 2 a^2\times \left(\mathfrak{Er}_1+\mathfrak{Er}_2\right),
\end{eqnarray}
where, for $i=1, 2$, $\displaystyle \mathfrak{Er}_i := \int_0^1 \E \Big [
\frac 1 {(n+1)^2}   \sum_{k,k'=0}^{n}
\big | c_{k'}-\E c_{k'}\big|
\times  \big ( 1\wedge\theta _{i,k,k'}^{-2}(\lambda)\big)
\Big]\,d\lambda$.
Thus $\ds \theta_{1,k,k'}(\lambda) =\delta_n' \big (
(k-k') - \lambda  k'z_n\big)$
with $\ds\delta_n' := \frac{\E \tau_n}{n} $ and $\ds z_n:=\frac{\tau_n-\E \tau_n}{\E \tau_n}$. Then, using $\delta'_n \limiten 0$, for $n$ large enough,
\ban
\mathfrak{Er}_1&=&
\int_0^1 \E \Big [
\frac C {(n+1)^{2}}   \sum_{k,k'=0}^{n}
\big | (k' \delta_n')\, z_n   \big|
\times  \Big ( 1\wedge\left[\delta_n' \big (
(k-k') - \lambda  k'z_n\big)\right]^{-2} \Big)
\Big]\,d\lambda\\
&\le&
\int_0^1\E \Big [
\frac C {(\E \tau_n)^{2}}    \int_0^{\E \tau_n}\int_0^{\E \tau_n}\hspace{-4mm} dxdy
\big | y\, z_n  \big|
\times  \Big ( 1\wedge\big[ (x-y) - \lambda y\, z_n \big]^{-2} \Big)
\Big]\,d\lambda.
\ean
But, for all $\lambda \in (0,1)$, one has
\ban
&&\hspace{-0.5cm} \frac 1 {(\E \tau_n)^{2}}    \int_0^{\E \tau_n}\int_0^{\E \tau_n}\hspace{-4mm}
 |y|
\times  \big ( 1\wedge\left[ (x-y) - \lambda y\, z_n \right]^{-2} \big)dx\,dy
\\
&&\hspace{2cm}=
\E \tau_n \int_0^{1}\int_0^{1} |v|
\times  \big ( 1\wedge(\E \tau_n)^{-2}\left[ (u-v) - \lambda v\, z_n \right]^{-2} \big)du\,dv
\\
&&\hspace{2cm}\leq 2\E \tau_n \int_0^{2}\int_0^{\infty} \big ( 1\wedge(\E \tau_n)^{-2}s^{-2} \big)\,ds\,dt \leq 4.
\ean
Therefore
$ \mathfrak{Er}_1 \le 4 \E |z_n|.$ Now, using the same method for $\mathfrak{Er}_2$, one obtains,
\ban
\E\tau_n \times \E \big [ (I'_n(a)-I_n(a))^2\big ]&\le& C \, \E \tau_n \times \E |z_n|  \\
&\le& C \big (\var (c_n-c_0)\big)^{1/2}\\
&\limiten & 0
\ean
from assumptions. This induces that CLT (\ref{TCL1}) holds.\\
~\\
Now the asymptotic expansion (\ref{SN}) can be proved. Consider first the deterministic case and since
$$
\sum_{k,k'=0}^{n}
e^{ i  (k-k') \alpha}=\Big |\frac {1-e^{i(n+1)\alpha }} {1-e^{i\alpha }}\Big |^2=\frac {\sin^2 ((n+1)\alpha /2)} {\sin^2(\alpha/2)}
$$
then,
\begin{eqnarray*}
S^2_n(a) & =& \frac {2\,a^2}{(n+1)^2}  \int_{\R^2} \big |
\widehat{\psi}(a \xi)\big |^2 f(\xi)\big |
\widehat{\psi}(a \xi')\big |^2 f(\xi') d\xi  d \xi'  \sum_{k,k'=0}^{n}
e^{ i  (k-k')  \frac {(c_n-c_0)} {n}  (\xi-\xi')}\\
& =&\frac {2\,a^2}{(n+1)^2}\int_{\R^2}  \big |
\widehat{\psi}(a \xi)\big |^2 f(\xi)\big |
\widehat{\psi}(a \xi')\big |^2 f(\xi')\frac {\sin^2\big (\frac {c_n-c_0} 2 \frac {n+1}n (\xi-\xi')\big )}
{\sin^2 \big (\frac {c_n-c_0}{2 n} \,(\xi-\xi')\big )} d\xi d \xi'  \\
& =&\frac {16 a^2}{c_n-c_0 }  \int_{\R_+^2} \hspace{-2mm}\big |
\widehat{\psi}(a z')\big |^2 f(z')\Big |
\widehat{\psi}\Big (a \big (z'+\frac {2 z}{c_n-c_0}\big ) \Big )\Big |^2 f\Big (z'+\frac {2 z}{c_n-c_0}\Big )
\frac {\sin^2 \big (\frac {n+1} n z \big )}
{n^2  \sin^2 \big (\frac {z}{n}\big )} dz d z'.
\end{eqnarray*}
Let us define $\displaystyle h_n(x):=\frac {\sin\big (\frac {n+1} n x \big )}
{n \sin \big (\frac {x}{n}\big )}$ and $\displaystyle h(x):=\frac {\sin x} x$. For all $(z,z')\in \R^2$,
$$
\Big |
\widehat{\psi}\Big (a \big (z'+\frac {2 z}{c_n-c_0}\big ) \Big )\Big |^2 f\Big (z'+\frac {2 z}{c_n-c_0}\Big ) \limiten
\big |\widehat{\psi}(a z')\big |^2 f(z')~~~\mbox{and}~~~ h^2_n(z)\limiten h^2(z).
$$
However Lebesgue Theorem cannot be applied. Denote $\nu(x):=|\psi(ax)|^2 f(x)$ for $x >0$.
From Assumptions F and W$(1,r)$ with $r>1/2$, $\nu$ is a differentiable function in $(0, \infty)$ and $\exists C>0$,  $ \forall z',x>0$, $|\nu'(z'+x)|\leq C \, |\nu'(z')|$. Then,
\ban
\Big |\nu \big (z'+\frac {2 z}{c_n-c_0} \big ) -\nu (z' ) \Big | & \leq &
\frac {2 z}{c_n-c_0} \, C\, |\nu'(z')|.
\ean
Moreover, $|h_n(z')| \leq 1$ for all $z'\in \R$, and
\ban
\int _{-n}^{n} h^2_n(z)\, dz& =& \frac {1}{n^2} \, \int _{-n}^{n} \sum_{k,k'=0}^{n}
e^{2i  (k-k') \frac {z} {n}}\, dz\\
& =&\frac {2}{n} \,  \sum_{0\leq k'< k \leq n} \, \frac {\sin \big (2 \, (k-k')\big )} {(k-k')} + \frac {1}{n^2}  \sum_{k=0}^{n} 2\,n\\
& =&\frac {2}{n} \,  \left\{\sum_{k=1}^{n} \,\big (n+1-k \big )\,  \frac {\sin \big (2 \, k\big )} {k}\right\} +2 \big ( \frac {n+1} n\big ).
\ean
Therefore,
\ba \label{limhn}
 \lim_{n\to \infty}\int _{-n}^{n} h^2_n(z)\, dz=2 \,  \left\{\sum_{k=1}^{\infty} \,  \frac {\sin \big (2 \, k\big )} {k}\right\}+ 2=\pi,
\ea
since
$ \frac {2}{n} \, \sum_{1\leq  k \leq n} k    \frac {\sin \big (2  k\big )} {k} \leq 4 \frac {\log n}{n} \limiten 0
$
and from Dirichlet Theorem, $ x-\pi=-2  \sum_{n \geq 1} \frac{\sin(n x)}n$ for all $x\in (0,2\pi)$.
Now, for $z'\geq 0$ and $n$ large enough,
\ban
\Big |\int_{\R_+} \upsilon \big (z'+\frac {2 z}{c_n-c_0} \big ) h^2_n(z) dz - \upsilon(z') \int _{0}^{n}  h^2_n(z)  dz\Big |&& \\
&& \hspace{-6.3cm} \leq \frac {2\, C \, |\upsilon(z')|} {c_n-c_0}  \int _{0}^{n}  z  h^2_n(z)  dz +  \int _{n}^\infty  \upsilon \big (z'+\frac {2 z}{c_n-c_0} \big ) h^2_n(z) dz \\
&& \hspace{-6.3cm} \leq\frac {2\, C \, |\upsilon(z')|} {c_n-c_0} \, \int _{0}^{n}  z \, \frac {4\sin^2(z)}{z^2} \, dz +
\int _{n}^\infty  \upsilon \big (z'+\frac {2 z}{c_n-c_0} \big ) \, dz\\
&& \hspace{-6.3cm} \leq    \frac {8 \,C \, |\upsilon(z')|} {c_n-c_0}\,   \Big ( \int _{0}^{1}  z  \, dz +\int _{1}^{n}  \frac {1}{z} \, dz \Big )+ C \, f(z')
\int _{n}^\infty  \Big |\psi\Big (a \big (z'+\frac {2 z}{c_n-c_0} \big )\Big ) \Big |^2  dz\\
\\
&& \hspace{-6.3cm} \leq    \frac {8 \,C \, |\upsilon(z')|} {c_n-c_0}\,   \Big ( \int _{0}^{1}  z  \, dz +\int _{1}^{n}  \frac {1}{z} \, dz \Big )+ C \, a^{-2r}\, f(z')
\int _{n}^\infty  \frac {C_\psi}{\big (1+ \big (z'+\frac {2 z}{c_n-c_0} \big )
\big )^{2r}} \, dz
\ean
under Assumption W$(1, 1, r)$. But when $r>1/2$,
$$ \int _{n}^\infty  \frac {1}{\big (1+ \big (z'+\frac {2 z}{c_n-c_0} \big )
\big )^{2r}} \, dz  \leq \frac {c_n-c_0} 2 \, \int _{2n/(c_n-c_0)}^\infty  \frac {1}{x  ^{2r}} \, dx \leq \frac {2^{2r-1}} {2(2r-1)} \, \frac {(c_n-c_0)^{2r}}{n^{2r-1}},
$$
and therefore there exists $C>0$ such that for all $z'>0$ and for $n$ large enough,
\ban
\Big |\int_{\R_+} \upsilon \big (z'+\frac {2 z}{c_n-c_0} \big ) h^2_n(z) dz - \upsilon(z') \int _{0}^{n}  h^2_n(z)  dz\Big |&& \\
&& \hspace{-6.3cm} \leq  C \,\Big (  \frac {1+\log (n)} {c_n-c_0}\, |\upsilon(z')|+  \,f(z')\, \Big ( \frac {c_n-c_0}{n^{1-\frac 1 {2r}}} \Big )^{2r} \Big )\\
&& \hspace{-6.3cm} \limiten 0
\ean
under assumptions of Proposition \ref{TCL:ponctuel2}. Finally, with (\ref{limhn}) in mind, one deduces that for all $z'\geq 0$,
$$
\lim_{n\to \infty} \int_{\R_+} \upsilon \big (z'+\frac {2 z}{c_n-c_0} \big )\, h^2_n(z) \, dz
=\lim_{n\to \infty} \upsilon(z')\, \int _{0}^{n}  h^2_n(z) \, dz=\frac{\pi} 2 \, \upsilon(z').
$$
Therefore, using the same method in $\R_-$ as in $\R_+$, one obtains
\begin{multline*}
\int_{\R}\big |
\widehat{\psi}(a z')\big |^2 f(z') \int_{\R} \, \Big |
\widehat{\psi}\Big (a \big (z'+\frac {2 z}{c_n-c_0}\big ) \Big )\Big |^2 f\Big (z'+\frac {2 z}{c_n-c_0}\Big )\,\frac {\sin^2 (z)}
{n^2 \, \sin^2 \big (\frac {z}{n}\big )}\, dz \, d z'\\
\limiten \pi
\, \int_{\R}\big |
\widehat{\psi}(a z')\big |^4 f^2(z')\, d z',
\end{multline*}
providing the asymptotic behavior of $S_n^2$. 
The proof is similar in the stochastic case with $c_n-c_0$ replaced by $\E (c_n-c_0)$.
\end{dem}
\vspace{-1cm}
\subsection*{Proof of Theorem \ref{Theo:TCL:discret}}
The proof of Theorem \ref{Theo:TCL:discret} uses the following
lemmas:
\begin{lem}\label{maj:cov}
Let $X$ be a Gaussian process defined by (\ref{repr:harmonizable})
with a spectral density $f$ satisfying (\ref{cond:f}) and
Assumption F($H$). Let us define, \ban \label{maj:cov:def} R(t,u,t',u'):=\E \Big [
(X(t+u)-X(t))\cdot (X(t'+u')-X(t'))~| ~{\cal F}_X \Big], \ean  for $(t,t')\in \R^2,\,(u,u')\in \R_+^2$. Then
$\exists C_f>0$ depending only on the spectral density $f$
such that for all $(u,u',t,t')\in
\R_+^2\times
\R^2$, with $\beta=\big(t'-t+u'-u\big)$,
\ban
\big|R(t,2u,t',2u')\big| \le C_f \Big  ( u \,
u'+(u \, u')^{(H+1)/2}|\beta|^{H-1}\1_{0<H<1}\big)
\Big).
\ean
\end{lem}
\begin{dem}
To begin with, remark that for  all $(t,t') \in \R^2$, $(u,u')\in \R^2_+$,
\begin{eqnarray*}
R(t,2u,t',2u')& =&\int _{\R}(e^{-i(t+2u)\xi}
-e^{-it\xi})(e^{i(t'+2u')\xi}
-e^{it'\xi})\,f(\xi) \, d\xi \nonumber \\
\nonumber & =& \int _{\R} (e^{-iu\xi} -e^{iu\xi})(e^{iu'\xi}
-e^{-iu'\xi})  \,e^{i\xi(t'-t)+i\xi(u'-u)} \,\,f(\xi)\,
d\xi \\
\nonumber & =& 8 \int _0 ^\infty  \sin (u\xi) \cdot \sin(u'\xi)
\cdot \cos \big ( \xi(t'-t+u'-u) \big )\,f(\xi) \, d\xi
\\
& =& 8 \left(I_1+I_2\right)
\end{eqnarray*}
with $ I_1:=\int _0 ^{\omega_K} \cdots d\xi$ and
$I_2=\int _{\omega_K}^\infty \cdots d\xi$. From the one
hand, with $|\sin a| \le |a|$ and $|\cos a| \le 1$, \ban |I_1|
&\le& u\, u' \, \int _0 ^{\omega_K} \xi^2\,f(\xi) \, d\xi \,\le\,
C\,u\, u', \ean where the last inequality follows from
(\ref{cond:f}). Now, if $H>1$, then the same bound can be extended to $I_2$ since $|\xi^2 f(\xi)|\leq C \xi ^{1-2H}$ and $\int_{\omega_K}^\infty \xi ^{1-2H}d\xi <\infty$. Therefore if $H>1$,
$$
|R(t,2u,t',2u')| \leq C_f\, u\, u'
$$
with $C_f$ only depending on $f$. \\
If $0<H<1$, firstly one obtains with a change of variable, $|\cos a| \le 1$ and $|\sin a|
\le (1 \wedge |a| )$
\ba
\nonumber |I_2| &\leq& \frac 1 {\sqrt{u\, u'}} \, \int_{\omega_K \sqrt{u\, u'}}^\infty |\sin(u\xi/\sqrt{u\, u'})\sin(u\xi/\sqrt{u\, u'})|\, f(\xi/\sqrt{u\, u'})d\xi \\
\nonumber  & \leq & (u\, u')^{H}  \, \int_{0}^\infty |\sin(u\xi/\sqrt{u\, u'})\sin(u\xi/\sqrt{u\, u'})|\, \xi^{-2H-1}d\xi \\
\nonumber  &\leq & (u\, u')^{H}  \Big ( \int_{0}^1 \xi^{-2H+1}d\xi + \int_{1}^\infty \xi^{-2H-1}d\xi \Big ) \\
\label{uu'2H}&\leq & C_H \, (u\, u')^{H},
\ea
with $C_H>0$ depending only on $H$. Secondly, with $\beta=\big(t'-t+u'-u\big)$ and an integration by parts,
\begin{multline*}
I_2 =  \int_{\omega_K}^\infty \sin \big( u \, \xi\big ) \sin  \big(u' \, \xi\big )\cos\big(\beta \, \xi\big
)f\big(\xi\big ) d\xi   \\
\hspace{1cm}=
\beta^{-1}\Big ( \Big[ \sin \big( u \, \xi\big ) \sin  \big(u' \, \xi\big )\sin\big(\beta \, \xi\big
)f\big(\xi\big ) \Big]_{\omega_K}^\infty
-\int_{\omega_K}^\infty \sin\big(\beta \, \xi\big
) \frac {\partial}{\partial x}\Big (\sin \big( u \, \xi\big ) \sin  \big(u' \, \xi\big ) \,f\big(\xi\big )\Big )dx \Big ),
\end{multline*}
where Assumption F($H$) insures the convergence of bracket term at
$\infty$. Using again
Assumption F($H$) for $f'$, changes of variables, $|\cos a| \le 1$ and $|\sin a|
\le (1 \wedge |a| )$,
\ban
I_2  & \leq &C_f \, u\, u'+
 C_0\, \beta^{-1}\int_{ \omega_K}^\infty \big |\sin\big
(\beta \, x\big )\big | \Big ( u\, \big |\sin
\big (u' \, x\big ) \big |+u' \big |\sin\big (u \, x\big )\big |\Big ) \,x^{-2H-1} dx\\
&& \hspace{5cm} +C'_0 \,\beta^{-1}\int_{ \omega_K}^\infty \Big |\sin\big
(\beta \, x\big ) \sin \big (u\, x\big )\sin
\big (u' \, x\big )\Big | x^{-2H-2} dx \\
 & \leq &C_f \, u\, u'+
 C_0\, \beta^{-1}(u \, u')^{H}\int_{0}^\infty \big |\sin (\beta \, x /\sqrt{u\, u'})\big | \Big ( u\, \big |\sin
\big (u' \, x/\sqrt{u\, u'} \big ) \big |+u' \big |\sin\big (u \, x/\sqrt{u\, u'} \big )\big |\Big ) \,x^{-2H-1} dx\\
&& \hspace{2cm} +C'_0 \,\beta^{-1}(u \, u')^{H+1/2}\int_{0}^\infty \big |\sin \big(\beta \, x /\sqrt{u\, u'} \big) \, \sin
\big (u' \, x/\sqrt{u\, u'}\big )\, \sin\big (u \, x/\sqrt{u\, u'} \big )\big |\Big ) \,x^{-2H-2} dx \\
 & \leq &C_f \, u\, u'+
C_0\big ( 2 \, (uu')^{H}\,  \int_{0}^1 x^{-2H+1} dx+ \beta^{-1}(u+u')(u \, u')^{H}\int_1^\infty x^{-2H-1} dx \big ) \\
&& \hspace{5cm} +C'_0 \big (  (u \, u')^H\,  \int_{0}^1 x^{-2H+1} dx +\beta^{-1}(u \, u')^{H+1/2}  \int_1^\infty x^{-2H-2} dx \big ) \\
& \leq & C_f \big ( u\, u'+(u \, u')^{H}+ \beta^{-1}(u+u')(uu')^{H}+ \beta^{-1}(u \, u')^{H+1/2}\big ),
\ean
where $C_f>0$ only depends on $f$. Therefore with (\ref{uu'2H}), for $0<H<1$,
$$|I_2| \leq C \, (u\, u')^{H} \wedge \big ( \beta^{-1}(u+u')(uu')^{H}+ \beta^{-1}(u \, u')^{H+1/2}\big ).$$
But the inequality
$(x\wedge y)\leq x^\alpha y^{1-\alpha}$ which is valid for all
$x,y
\geq 0$ and $0 \leq \alpha \leq 1$. Applied to previous inequality with appropriated choices of $\alpha$, one obtains $|I_2| \leq C \, (u\, u')^{(H+1)/2}\beta^{H-1}$. This
completes the proof of Lemma \ref{maj:cov}.
\end{dem}
Let $\varepsilon_n(a,b)$ be the error between the wavelet coefficient and its approximation, {\it i.e.}
\begin{eqnarray*}
\varepsilon_n(a,b): =  d_{X}(a,b)
- e^{}_{X}(a,b).
\end{eqnarray*}
The error $ \varepsilon_n(a,b)$ contains three different terms. The first one corresponds to the
replacement of the integral onto the interval $[0, T_n ]$ by its
Riemann sum. The second and the third ones correspond to the
replacement of the integral onto $\R$ by the integral onto the
interval $ [0, T_N ]$. More precisely, one has
\begin{eqnarray}\label{def:epsilon}
\varepsilon_n(a,b):=\varepsilon_{1,n}(a,b)+
\varepsilon_{2,n}(a,b) +\hspace{-0.5mm} \varepsilon_{3,n}(a,b),
\end{eqnarray}
with \ban
\varepsilon_{1,n}(a,b) &:=& a^{-1/2}\, \left\{\int_{ 0}^{T_n } \psi\Big
(\frac{t-b}{a} \Big)\,X(t)\,dt \,-\, \sum _{i=0}^{n-1}X
\big(t_i^{(n)}\big)\times \int_{t_i^{(n)}}^{t_{i+1}^{(n)}} \psi
\Big(
\frac{t-b}{a}\Big)\right\},
\\
\varepsilon_{2,n}(a,b) &:=&a^{-1/2}\,\int_{T_n  }^{\infty}
\psi\Big(\frac{t-b}{a} \Big)\,X(t)\,dt,
\\
\varepsilon_{3,n}(a,b) &:=& a^{-1/2}\,\int_{-\infty }^{ 0}
\psi\Big(\frac{t-b}{a} \Big)\,X(t)\,dt. \ean
The following lemma give bounds on $\E
\big | \varepsilon_{i,n}(a,k)\big |^2$ for $i=1,2, 3$.
\begin{lem}\label{maj:sum:epsilon}
Let $X$ be a Gaussian process defined by (\ref{repr:harmonizable})
with a spectral density $f$ satisfying (\ref{cond:f}) and
Assumption F($H$). Assume also Assumptions W$(1,3, 0)$ and let $(p_1,p_2)\in [1,\infty)^2$ and $(q_1,q_2)$ be defined by
$\frac 1 {p_j} + \frac 1 {q_j} = 1$ for $i=1, 2$.
 Then, there exist positive constants  $C_1$, $C_2$  and $C_3$ depending only on $f$, $a_{min}$, $a_{max}$, $C_\psi$ and $p_1,~p_2$, such that for
any r.v. $b$ independent on ${\cal
F}_X$ and satisfying  $\ds T_n^{\rho} \le
b\le T_n - T_n^{\rho}$ with $\rho > 1/2$, we have
\ba
\label{maj:E:eps1}  \E \big ( \big | \varepsilon_{1,n}(a,b)\big | ^2~| ~{\cal F}_X
\big) &\leq&  C_1 \,  \Big \{  \,\big\|\psi\big\|^2_{{\cal L}^{q_1}}
\Big(\sum_{i=0}^{n-1} L_i^{p_1+1} \Big)^{2/p_1} \, \delta_n^{2+2/p_1}
\\
\nonumber
&&
\nonumber
+ \1_{0<H<1}\cdot \|\psi\| _{{\cal L}^{q_2}}\big (\|\psi\| _{{\cal L}^{q_2}}+ \|\psi\|_{\infty} \big )
\Big (\sum_{i=0}^{n-1}  L_i^{1+
p_2 (1+H)/2} \Big )^{2/p_2}\delta_n^{1+H+2/p_2} \Big \};
 \\
\label{maj:E:eps23}
  \E \big (\big |  \varepsilon_{i,n}(a,b)\big |^2~| ~{\cal F}_X
\big)  &\le& C_2 \, 
a^{5}T_n^{2-4\rho}
~~~~\mbox{for $n$ large enough and  $i=2,3$}.
\ea
\end{lem}
\begin {rem}
When $H=1$, the term $\delta_n^{1+H+2/p_2}=\delta_n^{2+2/p_2}$ should be replaced by
$\ln\delta_n\times\delta_n^{2+2/p_2} =\mathcal{O}\left(\delta_n^{2-\varepsilon+2/p_2}\right)$ for all $\varepsilon>0$. Thus, if $H=1$, it suffices to replace $H$ by $1^-$ and formula (\ref{maj:E:eps1}) remains valid. This convention will be adopted in the following, in order to lighten the notations.
\end{rem}
\begin{dem} {\bf (1) Bound of $
\E \big (\big | \varepsilon_{1,n}(a,b)\big | ^2~| ~{\cal F}_X \big)$.}\quad
To begin with,
\begin{eqnarray*}
\E \big ( \big | \varepsilon_{1,n}(a,b)\big | ^2~|~ {\cal F}_X \big) && \\
&& \hspace{-3.2cm} = \frac 1 a  \sum_{i,j=0}^{n-1}
\int_{t^{(n)}_i}^{t^{(n)}_{i+1}}
\int_{t_j^{(n)}}^{t^{(n)}_{j+1}} \psi \Big
(\frac{t- b}{a} \Big  ) \overline{ \psi \Big (\frac{t'-b}{a}
\Big  )} \E
\Big (\big(X(t)-X\big( t^{(n)}_i\big)\big)\big(X(t')- X\big(t_j^{(n)}\big)\big)| {\cal F}_X  \Big ) dt dt'  \\
&& \hspace{-3.2cm} =\frac 1 a \sum_{i=0}^{n-1}
\sum_{j=0}^{n-1}\int_{t^{(n)}_i}^{t^{(n)}_{i+1}}
\int_{t_j^{(n)}}^{t^{(n)}_{j+1}} \psi \Big
(\frac{t-b}{a} \Big  ) \overline{  \psi \Big (\frac{t'-b}{a} \Big  )}
R\big(t^{(n)}_i,t-t^{(n)}_i,t_j^{(n)},t'-t_j^{(n)}\big)dt dt'.
\end{eqnarray*}
Lemma \ref{maj:cov}, with $2u= t-t^{(n)}_i$, $2u'= t'-t_j^{(n)}$ and $\beta=t^{(n)}_j-t^{(n)}_i+\frac 1 2 (t'-t^{(n)}_j)- \frac 1 2 (t-t^{(n)}_i)=\frac 1 2 (t'+t^{(n)}_j-t-t^{(n)}_i)$, implies
 \ban \E \big (\big | \varepsilon_{1,n}(a,b)\big | ^2~| ~{\cal F}_X \big)
&
\leq & C_f a^{-1} \big(S_1 +S_2 \1_{0<H<1}\big)\ean
with
\ban
S_1  & :=&  \sum_{i,j=0}^{n-1}
\int_{t^{(n)}_i}^{t^{(n)}_{i+1}}
\int_{t_j^{(n)}}^{t^{(n)}_{j+1}} u \, u' \Big |  \psi \Big
(\frac{t-b}{a} \Big  )  \Big |  \Big |  \psi \Big (\frac{t'-b}{a} \Big  ) \Big | dtdt' \\
S_2  & :=&  \sum_{i,j=0}^{n-1}
\int_{t^{(n)}_i}^{t^{(n)}_{i+1}}
\int_{t^{(n)}_{j}}^{t^{(n)}_{j+1}}  (u \, u')^{(H+1)/2}|\beta|^{H-1}
 \Big | \psi \Big
(\frac{t-b}{a} \Big  )  \Big |  \Big | \psi \Big (\frac{t'-b}{a} \Big  ) \Big |dtdt'. \ean
Set
\ba
\label{def:chi}
\chi(t) &:=& \sum_{i=0}^{n-1} |t-t^{(n)}_i| \1_{[t^{(n)}_i,t^{(n)}_{i+1}]}(t).
\ea
We have
$S_1 =
\left(\int_0^{T_n} \chi(t) \,\left|\psi \Big (\frac{t-b}{a}
\Big )\right| dt\right)^2$, then from H\"{o}lder Inequality, for $(p, q) \in
[1,\infty]^2$ with $1/p+1/q = 1$, we have
$$S_1 \le
\big\|\chi\big\|^2_{{\cal L}^p}\,
\big\|\psi \Big (\frac{t-b}{a}
\Big )\big\|^2_{{\cal L}^q}.$$
Obviously
$\ds \big\|\psi \Big (\frac{t-b}{a}
\Big )\big\|_{{\cal L}^q} = a^{1/q} \big\|\psi\big\|_{{\cal L}^q}$.
From the other hand, for $p<\infty$
\ban
\big\|\chi^H\big\|_{{\cal L}^p} & =& \Big ( \int_{\R} \chi^{Hp}(t)\, dt \Big )^{1/p}
\, =\,
\Big ( \int_{\R}  \sum _{i=0}^{n-1} \1_{[t^{(n)}_i,t^{(n)}_{i+1}]}(t)\,
 |t-t^{(n)}_i|^{Hp}
\, dt \Big )^{1/p}
\\
 & =&  \Big (\sum_{i=0}^{n-1} \int_{t^{(n)}_i}^{t^{(n)}_{i+1}} |t-t^{(n)}_i|^{Hp}
 dt \Big )^{1/p}
 \\
 &=& (1+Hp)^{-1/p}\,
\Big(\sum_{i=0}^{n-1} L_i^{1+Hp}\Big)^{1/p}\,\delta_n^{H+1/p}.
\ean
With the convention $1/\infty=0$, this result remains in force for $p=\infty$. It
follows for all $1 \leq p_1< \infty$,
\ba
\label{maj:ES1:stoch}
S_1 &\le& \frac {a^{2-2/p_1}\,\big\|\psi\big\|^2_{{\cal L}^{q_1}}}{
(p_1+1)^{2/p_1}}\,
\Big(\sum_{i=0}^{n-1} L_i^{p_1+1} \Big)^{2/p_1} \, \delta_n^{2+2/p_1}.
\ea
Next, in order to bound $S_2$ for $0<H <1$, write
\ban
S_{2}  & \leq &  C_f \, \sum_{0\leq i \leq j\leq
n-1}\int_{t^{(n)}_i}^{t^{(n)}_{i+1}}
\int_{t_j^{(n)}}^{t^{(n)}_{j+1}} \Big | \psi \Big
(\frac{t-b}{a} \Big  )  \psi \Big (\frac{t'-b}{a} \Big  ) \Big | (u
u')^{(1+H)/2}|\beta|^{-(1-H)}  dt  dt', \ean
where $u=\frac 1 2  (t-t^{(n)}_i)$, $u'=\frac 1 2  (t'-t^{(n)}_j)$ and $\beta =u'-u$. But since  $\beta=\frac 1 2 \big
( (t'+t_j^{(n)})-(t+t_{i})\big)$ for $i\leq j$ and $t \in
[t^{(n)}_i,t^{(n)}_{i+1}]$, $t' \in [t^{(n)}_{j},t^{(n)}_{j+1}]$ then $|\beta| \geq \frac 1 2 |t'-t|$. Therefore,,
\ban S_{2} &\leq & C_f   \sum_{0\leq
i\leq j\leq n-1}\int_{t^{(n)}_i}^{t^{(n)}_{i+1}}
\int_{t_j^{(n)}}^{t^{(n)}_{j+1}} \Big | \psi \Big
(\frac{t-b}{a} \Big  )   \psi \Big (\frac{t'-b}{a} \Big  ) \Big |(u
u')^{(1+H)/2} |t-t'|^{H-1}   dt\,dt'
\\
 &\leq &     \int_{0}^{T_n} \int_{0}^{T_n}
\chi(t)^{(1+H)/2}\chi(t')^{(1+H)/2}  \Big |
\psi \Big
(\frac{t-b}{a} \Big  )    \psi \Big (\frac{t'-b}{a} \Big  ) \Big |  |t-t'|^{H-1}  dt  dt'
\\
 &\leq &
\Big\|\chi(t)^{(1+H)/2} \chi(t')^{(1+H)/2}\Big\|_{{\cal L}^p(\R^2)}
\Big\| \psi \Big
(\frac{t-b}{a} \Big  )  \overline{  \psi \Big (\frac{t'-b}{a} \Big  )}  |t-t'|^{H-1}
\Big\|_{{\cal L}^q(\R^2)}
\ean
for any $(p, q) \in [1,\infty]^2$ with $1/p+1/q = 1$. But for all $p\ge 2$
 \ban \left\|\chi(t)^{(1+H)/2} \chi(t')^{(1+H)/2}\right\|_{{\cal L}^p(\R^2)}
 =
\Big(1+p(1+H)/2\Big)^{-2/p} \Big (\sum_{i=0}^{n-1} L_i^{1+ p (1+H)/2}\Big)^{2/p}
\ean
Next, with
$u=\big(t-c_k\big)/a$ and $v=\big(t'-c_k\big)/a$,
one gets
\ban
\Big\| \psi \Big
(\frac{t-b}{a} \Big  )    \psi \Big (\frac{t'-b}{a} \Big  )  |t-t'|^{H-1}
\Big\|_{{\cal L}^q(\R^2)}^q
&\le & a^{2+q(1-H)}
\int_{\R^2} \frac{|\psi(u)\psi(v)|^q}{|u-v|^{q(1-H)}} dudv
\\
&\hspace{-13cm}\le&
\hspace{-6.8cm}a^{2+q(1-H)}\Big (\int_{\R^2,|u-v|\ge 1} |\psi(u)\psi(v)|^q dudv
+\int_{\R^2, |u-v|< 1}
\frac{|\psi(u)\psi(v)|^q}{|u-v|^{q(1-H)}} dudv\Big)
\\
&\hspace{-13cm}\le&
\hspace{-6.8cm}a^{2+q(1-H)}\Big ( \|\psi\|_{{\cal L}^q}^{2q} + \, \|\psi\|_{\infty}^{q}
\,\|\psi\|_{{\cal L}^q}^{q}\int_0^1
s^{-q(1-H)}ds\Big).
\ean
The last integral is equal to $\big(1-q(1-H)\big)^{-1}$ when $p>1/H$.
Thus, for all $ p_2>1/H$,
\ba
\label{maj:ES4:stoch}
S_{2}
&\le &  \frac {\|\psi\|_{{\cal L}^{q_2}}
\big (\|\psi\|_{{\cal L}^{q_2}}^{q_2}  + \frac {\|\psi\|_{\infty}^{q_2}}{1-q_2(1-H)}
\big )^{1/q_2}}{a^{H-3+2/p_2} \big(1 + p_2(1 + H)/2\big)^{2/p_2}}
\Big (\sum_{i=0}^{n-1}  L_i^{1+
p_2 (1+H)/2} \Big )^{2/p_2}\delta_n^{1+H+2/p_2}.
\ea
Finally by summing up (\ref{maj:ES1:stoch}) and
(\ref{maj:ES4:stoch}), one gets the bounds of $\E \big (\big |  \varepsilon_{1,n}(a,b)\big |^2~| ~{\cal F}_X \big)$. \\
~\\
{\bf (2) Bound of $\ds \E \big (
\big | \varepsilon_{2,n}(a,b)\big | ^2~| ~{\cal F}_X \big)$.} Since $T_n$ is  independent on ${\cal F}_X$,
 \ban \E \big ( \big | \varepsilon_{2,n}(a,b)\big | ^2~| ~{\cal F}_X \big) &=& \frac 1 a \int_{T_n
}^{\infty}\int_{T_n }^{\infty}
\psi \Big
(\frac{t-b}{a} \Big  )  \overline{  \psi \Big (\frac{t'-b}{a} \Big  )}\E \big (
X(t)X(t')\big )dtdt'
\\
&\le&  \frac {C_f} a \Big (\int_{T_n }^{\infty}
\left|\psi\Big(\frac{t-b}{a} \Big)\right|\big(1+|t|\big)dt
\Big)^2.\ean
from Lemma
\ref{maj:EXt:grand}. But, according to Assumption W$(1,3)$, $(1+|t|^3) |\psi(t)|$ is a bounded function and
\ban \E \big ( \varepsilon_{2,n}^2(a,b)~| ~{\cal F}_X \big)
& \le& C_f\, C_\psi^2\, a^{-1} \,\Big(\int_{T_n }^{\infty}
 \big(1+t\big) \big(1+ \big(t-b\big)/a\big )^{-3}
dt\Big )^2. \ean
If $T_n\ge 1$, then $1+ \big(t-b\big)/a \le 1+
\big(t-T_n + T_n^\rho\big)/a$ for all $t\ge T_n$ and with the change of variable $\ds v =
\big(t-\,T_n+T_n^\rho\big)/T_n$,
\ban \int_{T_n }^{\infty} \frac {\big(1+t\big)}{\big(1+
\big(t-b\big)/a \big)^{3}}\,dt &\le& T_n\, a^3\,
\int_{T_n^{\rho-1}}^{\infty}  \frac {\big(1+  \,v\,T_n+ \,T_n-T_n^\rho \big)}{\big(a+
 \,v\,T_n \big)^{3}}\,dv
\\
&\le& T_n^{-1}\, a^3\, \int_{T_n^{\rho-1}}^{\infty}  \frac
{\big(v+2
\big)}{v^{3}}\,dv
\;=\,
a^3\,\big[ T_n^{1-2\rho}+ T_n^{-\rho}\big]
\ean
If $T_n\le 1$, by using $b \le T_n$,
\ban \int_{T_n }^{\infty} \frac {\big(1+t\big)}{\big(1+
\big(t-b\big)/a \big)^{3}}\,dt
\;\le\;
\int_{0}^{\infty}
\frac {\big(1+v + T_n\big)}{\big(1+
v/a \big)^{3}}\,dv
\leq  a + \frac 1 2 \,a^2.
\ean
Eventually, one deduces the bound of $\ds \E \big (
\big | \varepsilon_{2,n}(a,b)\big | ^2~| ~{\cal F}_X \big)$.\\
\\
{\bf (3) Bound of $\ds \E \big (
\big | \varepsilon_{3,n}(a,b)\big | ^2~| ~{\cal F}_X \big)$.}
Find a bound for $\ds \E \big (
\big | \varepsilon_{3,n}(a,b)\big | ^2~| ~{\cal F}_X \big)$
follows the same steps than for bounding $\ds \E \big (
\big | \varepsilon_{2,n}(a,b)\big | ^2~| ~{\cal F}_X \big)$.
\end{dem}
\begin{lem}
\label{lem:maj:vn:stoch}
Under assumptions of Lemma \ref{maj:sum:epsilon} and if $\ds s> 2 +\frac 1 {2H} \, \big [1-3H\big ]_+ $ and if
\ban
\delta_n=n^{-d}\quad\mbox{with}\quad \Big ( \frac {1+(2H\wedge 1)}{1+s(2H\wedge 1)}\big )\vee \big (\frac {s+(H\wedge 1)}{s(2+(H\wedge 1))-1}  \Big )<d<1,
\ean
then for all $a>0$,
\ba\label{cond_vn}
(n\, \delta_n) \, \E\big ( \big |
\varepsilon_{n}(a,c_k)\big |^2 \big) \limiten 0\quad\mbox{for all $0 \leq k \leq n$}.
\ea
\end{lem}
\begin{dem}
With
$(x+y+z)^2
\le 3\,(x^2+y^2+z^2)$ for all real numbers $x$, $y$, $z$,
$$\ds
\E\big ( \big |
\varepsilon_{n}(a,c_k)\big |^2 \big)\,\le\, 3\,  \E \big [\E \big (\big |  \varepsilon_{1,n}(a,b)\big |^2~| ~{\cal F}_X
\big) \big ] + 6\, \E \big [\E \big (\big |  \varepsilon_{2,n}(a,b)\big |^2~| ~{\cal F}_X
\big) \big ].
$$
Then using
Lemma \ref{maj:sum:epsilon}, with $p_1\geq 1$ and $p_2>1/H$, an optimal choice of $p_1, p_2$ depends on $s$.  Hence, since $\E \left(|Z|^\alpha\right)\le \left(\E|Z|\right)^\alpha$
for any r.v. $Z$ and $\alpha \in [0,1]$ from Jensen Inequality,
\begin{eqnarray*}
&&\hspace{-5mm}\mbox{1. if $3\leq s$},~~\mbox{with}~~1+p_1=s, \\
&& \hspace{1cm} \E\Big(\sum_{i=0}^{n-1} L_i^{p_1+1} \Big)^{\frac 2 {p_1}} \delta_n^{2+\frac 2 {p_1}} \leq M_s^{\frac 2 {s-1}} \cdot (n \, \delta^s_n)^{\frac 2 {s-1}};\\
&&\hspace{-5mm}\mbox{2. if $\max \big (2+H,\frac 3 2 + \frac 1 {2H} \big) \leq s$},~~\mbox{with}~~1+\frac 1 2 \, p_2(1+H)=s, \\
&& \hspace{1cm} \E\Big(\sum_{i=0}^{n-1} L_i^{1+
p_2 (1+H)/2} \Big)^{\frac 2 {p_2}} \delta_n^{1+H+\frac 2{p_2}}  \leq M_s^{\frac {1+H} {s-1}} \cdot (n \,  \delta_n^s )^{\frac {1+H} {s-1}}.
\end{eqnarray*}
However, these inequalities may be extended to smaller values of $s$ by using the sharper inequality
$\ds \E \big (\sum |x_i|\big )^{\alpha \beta} \leq \E \big ( \sum |x_i|^\beta \big) ^\alpha \leq n^\alpha \big (\max _{0\leq i \leq n-1} \E (|x_i|^\beta) \big )^\alpha$ when $(\alpha,\beta)\in (0,1]^2$ and therefore for $r>1$, $\ds \E \big (\sum L_i^r\big )^{\alpha \beta}  \leq n^\alpha \big (\max _{0\leq i \leq n-1} \E (|L_i|^{r\beta}) \big )^\alpha$; it is then possible to choose $r\beta=s$. Thus,
\begin{eqnarray*}
&&\hspace{-5mm}\mbox{1'. if $2<s\leq 3$},~~\mbox{with~~$\alpha\beta=\frac 2 {p_1}$ and $s=\beta+\frac 2 \alpha$, the best possible choice is $\alpha=1$ and $\beta=s-2$}, \\
&& \hspace{1cm} \E\Big(\sum_{i=0}^{n-1} L_i^{1+\frac 2 {\alpha \beta}} \Big)^{\alpha \beta} \delta_n^{2+\alpha \beta} \leq M_s ^s  \cdot n \delta_n^s;\\
&&\hspace{-5mm} \mbox{2'. if $H\geq 1/2$, and $1+H<s\leq 2+H$},~~\mbox{with~~$\alpha\beta=\frac 2 {p_2}$ and $s=\beta+\frac{1+H}\alpha$, the best possible choice is} \\
&& \hspace{10cm}  \mbox{ $\alpha=1$ and $\beta=s-(H+1)$},  \\
&& \hspace{1cm} \E\Big(\sum_{i=0}^{n-1} L_i^{1+ \frac {1+H}{\alpha \beta} } \Big)^{\alpha \beta} \delta_n^{1+H+\alpha \beta} \leq M_s ^{s} \cdot  n \delta_n^s;\\
&&\hspace{-5mm} \mbox{2''. if $0<H\leq 1/2$ and  $\frac 1 2 + \frac {1}{ 2H}<s\leq \frac 3 2 + \frac {1}{ 2H}$},~~\mbox{with~~$\alpha\beta=\frac 2 {p_2}$ and $s=\beta+\frac{1+H}\alpha$, the best possible choice is}\\ && \hspace{10cm}  \mbox{ $\alpha=2H$ and $\beta=s-\frac {H+1} {2H}$,}  \\
&& \hspace{1cm} \E\Big(\sum_{i=0}^{n-1} L_i^{1+ \frac {1+H}{\alpha \beta} } \Big)^{\alpha \beta} \delta_n^{1+H+\alpha \beta} \leq M_s ^{2H}  \cdot  \big (n \delta_n^s\big )^{2H}.
\end{eqnarray*}
We finally obtain for $n$ large enough and using $n \, \delta_n \limiten \infty$ and $n \, \delta^{s}_n \limiten 0$,
\begin{itemize}
\item if $H\geq 1/2$, $\E \ds \big [ \E \big (
\big | \varepsilon_{1,n}(a,b)\big | ^2~| ~{\cal F}_X \big)\big ] \leq C \,\big (  n \delta_n^{s}\,  \1 _{2<s \leq 2+(H\wedge1)} +(n \delta_n^{s} ) ^{\frac {1+(H\wedge1)} {s-1}}\, \1 _{2+(H\wedge1) \leq s} \big )$;
\item if $0< H\leq  1/2$, $\E \ds \big [ \E \big (
\big | \varepsilon_{1,n}(a,b)\big | ^2~| ~{\cal F}_X \big)\big ] \leq C \,\big (  n \delta_n^{s}\,  \1 _{2\vee (\frac 1 2 + \frac 1 {2H})<s \leq \frac 3 2 + \frac 1 {2H}} +(n \delta_n^{s} ) ^{\frac {1+H} {s-1}}\, \1 _{\frac 3 2 + \frac 1 {2H} \leq s} \big )$;
\end{itemize}
Both these  inequalities may be reduced to only one for all $s > 2\vee (\frac 1 2 + \frac 1 {2H})$ and $H>0$:
\ban
\E \ds \big [ \E \big (
\big | \varepsilon_{1,n}(a,b)\big | ^2~| ~{\cal F}_X \big)\big ] \leq \left \{ \begin{array}{ll}
 C \,  \big (n \delta_n^{s} \big) ^{ 2H\wedge 1 } &\mbox{if  $2\vee (\frac 1 2 + \frac 1 {2H})< s<(2+H\wedge 1) \vee (\frac 3 2 + \frac 1 {2H})$}
\\
C \, \big (n \delta_n^{s} \big) ^{\frac {1+(H\wedge1)} {s-1} } & \mbox{if  $(2+H\wedge 1) \vee (\frac 3 2 + \frac 1 {2H}) \leq s$}
 \end{array}\right . .
\ean
Hence, with $\delta_n=C_\delta \, n^{-d}$,
\ban
&\bullet& \mbox{$\frac {1+(2H\wedge 1)}{1+s(2H\wedge 1)}\leq d<1$} \hspace{2cm}\mbox{  if $2 \vee \big (\frac 1 2 +\frac 1 {2H}\big )  \leq   s   < \big (2+(H\wedge 1)\big ) \vee \big (\frac 3 2 +\frac 1 {2H}\big ) $} \\
&\bullet&\mbox{$\frac {s+(H\wedge 1)}{s(2+(H\wedge 1))-1}\leq d<1$} \hspace{1.55cm} \mbox{  if $\big (2+(H\wedge 1)\big ) \vee \big (\frac 3 2 +\frac 1 {2H}\big )  \leq   s  \leq  \infty$}
\ean
then
\ban
\label{inegfinale}
(n \delta_n) \E \ds \big [ \E \big (
\big | \varepsilon_{1,n}(a,b)\big | ^2~| ~{\cal F}_X \big)\big ] \limiten 0.
\ean
To finish the proof of Lemma \ref{lem:maj:vn:stoch} it remains to
show $\ds (n  \delta_n)\E \ds \big [ \E \big (
\big | \varepsilon_{2,n}(a,b)\big | ^2~| ~{\cal F}_X \big)\big ]  \limiten 0$.\\
From Lemma \ref{maj:sum:epsilon} it follows that $\ds \E \ds \big [ \E \big (
\big | \varepsilon_{2,n}(a,b)\big | ^2~| ~{\cal F}_X \big)\big ] \le C \, \int_0^\infty g(x)\,f_n(x)\,dx$ where
$f_n$ is the probability distribution function of $T_n$
and $\ds g(x) =\1_{(x< 1)}+
\1_{(x\ge 1)} \,x^{2-4\rho} $. Since $\rho>3/4$, $g(x) \le 1$ for all
$x>0$ and $g$ is a non increasing map,
\ban
\int_0^\infty g(x)\,f_n(x)\,dx
&\le&
\int_0^{\frac 1 2 \,m_s\, n \,\delta_n}f_n(x)\,dx+g\big(\frac 1 2 \,m_s\, n \, \delta_n\big)\,
\int_{\frac 1 2 \,m_s \, n \,\delta_n}^\infty f_n(x)\,dx
\\
&\le& \Pr\big(T_n\le \frac 1 2 \,m_s\, n \, \delta_n\big) + g\big(\frac 1 2 \,m_s \, n \,
\delta_n\big)
\\& \leq & \Pr\big(|T_n-\E [T_n]|\geq\E [T_n]- \frac 1 2 \,m_s\, n \, \delta_n\big) + \big(\frac 1 2 \,m_s\, n \,
\delta_n\big)^{2-4\rho},\\
& \leq & 4 \, \frac {M_s} {m_s} \, \frac {n \, \delta^2_n} {n^2 \, \delta^2_n}
 + \big(\frac 1 2 \,m_s\, n \, \delta_n\big)^{2-4\rho},
\ean
from Bienaym\' e-Chebyshev Inequality since $s\geq 2$ and $\var (T_n)\leq M_s  n \delta_n^2$ from the independence of
$(L_i)_{i\in \N}$. Therefore $\ds
( n \, \delta_n)\,\E \ds \big [ \E \big (
\big | \varepsilon_{2,n}(a,b)\big | ^2~| ~{\cal F}_X \big)\big ] \limiten 0$.
\end{dem}
\begin{dem}[Theorem \ref{Theo:TCL:discret}]
Denote $\;\ds v_n(a)= \big (n \delta_n^{s} \big) ^{ 2H\wedge 1 }+ \big (n \delta_n^{s} \big) ^{\frac {1+(H\wedge1)} {s-1} }
+ \big(n\, \delta_n\big)^{2-4\rho}$.
Then, following the same method that in 
Bardet and Bertrand (2007b), pp. 33-35, one obtains
\ba \E\left|I_n(a) - J_n(a) \right| \le C \, v_n(a)^{1/2},\ea
and from this, Lemmas \ref{lem:maj:vn:stoch} and Slutsky Lemma,
the proof is achieved.
\end{dem}
\vspace{-1cm}
\subsection*{Proof of Proposition \ref{cvg:Ip}}
\begin{dem}
It is obvious that
\ban
\mathcal{I}_\lambda(a)\hspace{-2.6mm} &=& \hspace{-2.6mm}  \int_{\R} \big |\widehat{\psi}_\lambda(\xi)\big|^2
f(\xi/a) d\xi
= \lambda \int_{\R} \big|\widehat{\psi}\big(\lambda(\xi-1)\big)\big|^2
f(\xi/a) d\xi
\\
\hspace{-2.6mm}&=&\hspace{-2.6mm} \int_{\R} \big|\widehat{\psi}\big(v\big)\big|^2
f\Big (\frac 1 a + \frac v {a  \lambda}\Big ) dv.
\ean
Then, from a usual Taylor expansion, and since $\widehat \psi$ is supposed to be an even function supported in $[-\Lambda,\Lambda]$,
\ba
\nonumber \Big |\mathcal{I}_\lambda(a)- \big\|\widehat{\psi}\big \|^2_{{\cal L}^2(\R)} \, f(1/a)\Big | \hspace{-2.6mm}&\le&\hspace{-2.6mm} \frac 1 {2a^2 \lambda^2} \,
\Big (\sup_{-\Lambda /\lambda \leq h }\Big\{\Big|f''\big (\frac {1+h} a\big )\Big| \Big \}   \int_{-\Lambda}^{\Lambda}v^2
\big|\widehat{\psi}\big(v\big)\big|^2 dv \Big ).
\ea
For $\lambda > 2 \Lambda$, then $\ds \sup_{-\Lambda /\lambda \leq h }\Big\{\Big|f''\big (\frac {1+h} a\big )\Big| \Big \}\leq \sup_{x>1/2a}\big\{\big|f'' (x)\big| \big \}<\infty$. Therefore, since $\psi$ satisfies Assumption W$(1,5)$, there exists $C>0$ such that,
\ba
\label{majIp} \Big |\mathcal{I}_\lambda(a)- \big\|\widehat{\psi}\big \|^2_{{\cal L}^2(\R)} \, f(1/a)\Big |
 \hspace{-2.6mm}&\le&\hspace{-2.6mm} C \,  \frac 1 { \lambda^2}.
\ea
Let us denote $I^{(\lambda)}_n(a)$ (respectively $\mathcal{I}_\lambda(a)$, $\beta^{(\lambda)}_n$ and $S^{(\lambda)}_n(a)$) instead on $I_n(a)$ (resp. $\mathcal{I}_1(a)$, $\beta^{}_n$ and $S^{}_n(a)$)
when $\psi$ is replaced by $\psi_{\lambda}$. Firstly,
\begin{eqnarray*}
\frac 1 \lambda \,\Big (4\pi \, a^2 \, \int_{\R}\big |
\widehat{\psi_\lambda}(a z)\big |^4 f^2(z)\, d z \Big ) \hspace{-2.6mm}&=&\hspace{-2.6mm} 4\pi \,a \, \int_{\R}\big |
\widehat{\psi}(u)\big |^4 f^2\Big (\frac 1 a + \frac u{a\lambda} \Big )\, d u\\
\hspace{-2.6mm}& \limitep & \hspace{-2.6mm}4\pi \,a \,f^2(1/a) \, \int_{\R}\big |
\widehat{\psi}(u)\big |^4 \, d u,
\end{eqnarray*}
from Lebesgue Theorem. Hence, if $(\lambda_n)$ is a sequence such that
$\lambda_n \limiten \infty$,
\begin{eqnarray}\label{Sp} \frac{ \E (T_n)} {\lambda_n} \,  \big (S^{(\lambda_n)}_n(a)\big )^2 \limiten 4\pi \,a \,f^2(1/a) \, \int_{\R}\big |
\widehat{\psi}(u)\big |^4 \, d u.
\end{eqnarray}
Secondly, from the proof of Proposition \ref{TCL:ponctuel} and inequalities (\ref{bound_betan}) and (\ref{bound_Sn}),
there exists $C>0$ not depending on $n$ and $\lambda$,
$$
\beta^{(\lambda)}_n/S^{(\lambda)}_n \leq C \, \mathcal{I}^{-1}_\lambda(a) \, \big(n\,\max_{1 \leq k\leq n} \big|c_{k+1}
-c_{k}\big|\big)^{1/2-1/q}\quad\mbox{for all $q \in (1,2)$}.
$$
Thus, since $\mathcal{I}_\lambda(a)$ is bounded, $\beta^{(\lambda_n)}_n/S^{(\lambda_n)}_n  \limiten  0$ and Proposition \ref{TCL:ponctuel}
becomes:
$$
\frac{I^{(\lambda_n)}_n(a)
- \mathcal{I}_{\lambda_n}(a)}{S^{(\lambda_n)}_n(a)}\limiteloin \mathcal{N}(0,1).
$$
Finally, using (\ref{majIp}) and (\ref{Sp}), on deduces that for all $a>0$,
\ban \label{Ipp}
\sqrt {\frac{ \E T_n} {\lambda_n} } \big (I^{(\lambda_n)}_n(a)
- \big\|\psi\big \|^2_{{\cal L}^2(\R)}  f(1/a) \big ) \limiteloin  \mathcal{N}\Big (0 ,4\pi a f^2(1/a)   \int_{\R} \big |
\widehat{\psi}(u)\big |^4  d u\Big ),
\ean
when $(\lambda_n)_n$ is such that $\sqrt {\frac{ \E T_n} {\lambda_n} }\,\frac 1 {\lambda^2_n} \limiten 0$, {\it i.e.}
when $\ds \lambda_n^{-5} \, n \, \delta_n \limiten 0$. Since also $\ds \lambda_n^{-1}  n  \delta_n \limiten \infty$ (to obtain a consistent estimator), then with $\delta_n=n^{-d}$ and $\lambda_n=n^{d'}$,
\ba \label{cond1}
\frac {1-d}5 <d'<1-d.
\ea
Moreover, Proposition \ref{TCL:ponctuel2} has also to be checked. In its proof, $\E \tau_n$ has to be replaced by $\E \tau_n/\lambda_n$ and since the bounds $C \, (1\wedge|\theta|^{-1})$ in Lemma \ref{lem:gamma:a} have to be replaced by $C/\lambda_n^2
\, (1\wedge|\theta|^{-1})$, then condition $\ds n \delta_n^2 \limiten 0$ has to be replaced by  $\ds n \delta_n^2/\lambda_n^5 \limiten 0$, that is $\ds d'>\frac {1-2d}5$ which is satisfied when (\ref{cond1}) is satisfied.\\
~\\
It remains to control $\varepsilon^2_n(a,c_k)$ with Lemma \ref{maj:sum:epsilon} and \ref{lem:maj:vn:stoch}.
For all $1\leq q \leq \infty $, with $1/\infty=0$ by convention,
$$
\big \|{\psi}_\lambda\big \|_{{\cal L}^q}=\lambda^{(2-q)/2q} \,
\big \|{\psi}\big \|_{{\cal L}^q}\quad\mbox{and}\quad \big \|\widehat{\psi}_\lambda\big \|_{{\cal L}^q}=\lambda^{(q-2)/2q} \,
\big \|\widehat{\psi}\big \|_{{\cal L}^q}.
$$
Then, Lemma \ref{maj:sum:epsilon} becomes (with $\lambda_n \limiten \infty$):
\ban
\E \big ( \big | \varepsilon_{1,n}(a,b)\big | ^2~| ~{\cal F}_X
\big) \leq \left\{ \begin{array}{ll}
\displaystyle \bullet~ C_1  \,\|\psi\|^2_{{\cal L}^1}\,  \lambda_n \,  \delta_n^{1+(H\wedge 1)} &  \mbox{if $s=\infty$};\\
 \displaystyle \bullet~C_1   \Big \{  \|\psi\|^2 _{{\cal L}^{q_2}} \,  \lambda_n^{\frac 2 {q_2}-1} \,
\Big (\sum_{i=0}^{n-1}  L_i^{1+
p_2 (1+H)/2} \Big )^{2/p_2}\, \delta_n^{1+H+2/p_2} \1_{0<H<1}  & \\
\displaystyle \hspace{3.5cm}+\big\|\psi\big\|^2_{{\cal L}^{q_1}}\,  \lambda_n^{\frac 2 {q_1}-1} \,
\Big(\sum_{i=0}^{n-1} L_i^{p_1+1} \Big)^{2/p_1} \delta_n^{2+2/p_1}
 \Big \},&  \mbox{if $s<\infty$.}
 \end{array}
 \right .
\ean
The case $s<\infty$ can be more detailed following the values of $p_1$ and $p_2$ considered in Lemma \ref{lem:maj:vn:stoch}:
\ban
&&\hspace{-6mm} a/~\mbox{if $H\geq 1/2$ and} \\
&&\hspace{+6mm} \bullet ~\mbox{$2 < s \leq 2+(H\wedge 1)$, then $\frac  2 {p_1} =s-2$ and $\frac  1 {p_2}=(s-(H\wedge 1)-1)/2$ and therefore}\\
&& \hspace{6mm}\E \big (
\big | \varepsilon_{1,n}(a,b)\big | ^2~| ~{\cal F}_X \big) \le C  \big ( \lambda_n^{3-s}+\lambda_n^{2+(H\wedge 1)-s}\big ) \big (n \, \delta_n^s\big ) \\
&& \hspace{5.3cm} \Longrightarrow ~~ \E \big (
\big | \varepsilon_{1,n}(a,b)\big | ^2~| ~{\cal F}_X \big) \le C' \,  \lambda_n \,  \big (n \, \delta_n^s\big );\\
&&\hspace{+6mm} \bullet ~\mbox{$ 2+(H\wedge 1)\leq s$, then $p_1=s-1$ or $\frac  2 {p_1}=s-2$ and $\frac  1 {p_2}=(1+(H\wedge 1))/(2(s+1))$ and therefore}\\
&& \hspace{6mm}\E \big (
\big | \varepsilon_{1,n}(a,b)\big | ^2~| ~{\cal F}_X \big) \le C  \Big (\lambda_n^{3-s}\big (n \, \delta_n^s\big )+  \lambda_n^{\frac {s-3}{s-1}}\big (n \, \delta_n^s\big )^{\frac 2 {s-1}}  +\lambda_n^{\frac {s-(2+(H\wedge 1))}{s-1}}\big (n \, \delta_n^s\big ) ^{\frac {1+(H\wedge 1)} {s-1}}   \Big )\\
&& \hspace{5.3cm} \Longrightarrow ~~ \E \big (
\big | \varepsilon_{1,n}(a,b)\big | ^2~| ~{\cal F}_X \big) \le C' \,   \lambda_n\,  \big (n \, \delta_n^s\big ) ^{\frac {1+(H\wedge 1)} {s-1}}.\\
&&\hspace{-6mm} b/~\mbox{if $0<H\leq 1/2$ and} \\
&&\hspace{+6mm} \bullet ~\mbox{$2\vee \big (\frac 1 2 + \frac 1 {2H}\big )  < s \leq  \frac 3 2+\frac 1 {2H}$, then $\frac  2 {p_1}=s-2$ or $\frac 2 {p_1}=\frac 2 {s-1}$ and $\frac  1 {p_2}=H\big (s-(\frac 1 2 +\frac 1 {2H})\big )$ and therefore}\\
&& \hspace{6mm}\E \big (
\big | \varepsilon_{1,n}(a,b)\big | ^2~| ~{\cal F}_X \big) \le C  \big ( \lambda_n^{3-s}\big (n \, \delta_n^s\big ) + \lambda_n^{\frac {s-3}{s-1}}\big (n \, \delta_n^s\big )^{\frac 2 {s-1}} +\lambda_n^{1-2H(s-(\frac 1 2 +\frac 1 {2H}))}\big (n \, \delta_n^s\big )^{2H}\big ) \\
&& \hspace{5.3cm} \Longrightarrow ~~ \E \big (
\big | \varepsilon_{1,n}(a,b)\big | ^2~| ~{\cal F}_X \big) \le C' \,  \lambda_n \big (n \, \delta_n^s\big )^{2H};\\
&&\hspace{+6mm} \bullet ~\mbox{$ \frac 3 2+\frac 1 {2H} \leq s $, then $\frac  2 {p_1}=s-2$ or $\frac  2 {p_1}=\frac 2 {s-1}$ and $\frac  1 {p_2}=(1+H)/(2(s+1))$ and therefore}\\
&& \hspace{6mm}\E \big (
\big | \varepsilon_{1,n}(a,b)\big | ^2~| ~{\cal F}_X \big) \le C  \big ( \lambda_n^{3-s}\big (n \, \delta_n^s\big )+\lambda_n^{\frac {s-3}{s-1}}\big (n \, \delta_n^s\big )^{\frac 2 {s-1}}  + \lambda_n^{1-2H(s-(\frac 1 2 +\frac 1 {2H}))}\big (n \, \delta_n^s\big )^{2H}\big ) \\
&& \hspace{5.3cm} \Longrightarrow ~~ \E \big (
\big | \varepsilon_{1,n}(a,b)\big | ^2~| ~{\cal F}_X \big) \le C' \,  \lambda_n \big (n \, \delta_n^s\big )^{\frac {1+H} {s-1}}.
\ean
Note that the bound is not always optimal but it simplifies a lot of different subcases. Condition (\ref{cond_vn}) is now
$\ds \frac {n \, \delta_n} {\lambda_n}\, \E \big (
\big | \varepsilon_{1,n}(a,b)\big | ^2~| ~{\cal F}_X \big) \limiten 0$. Therefore in any case this condition does not depend on $\lambda_n$ and it can be summarize with $\delta_n=n^{-d}$ with (the case $s=\infty$ is obtained by replacing with the limit):
\ba
\label{cond3}\bullet&~~ \mbox{if $2 \vee \big (\frac 1 2 +\frac 1 {2H}\big )  \leq   s   < \big (2+(H\wedge 1)\big ) \vee \big (\frac 3 2 +\frac 1 {2H}\big ) $,~~ }&d > \frac {1+(2H\wedge 1)}{1+s(2H\wedge 1)};\\
\label{cond31}\bullet& \hspace{-1.6cm} \mbox{if $\big (2+(H\wedge 1)\big ) \vee \big (\frac 3 2 +\frac 1 {2H}\big )  \leq   s  <  \infty$,~~}&d > \frac {s+(H\wedge 1)}{s(2+(H\wedge 1))-1}.
\ea
Finally, for $b \leq  T_n-T_n^\rho$, with $\psi$ satisfying Assumption W$(1,3,1)$:
\ban \E \big ( \varepsilon_{2,n}^{\lambda,2}(a,b)~| ~{\cal F}_X \big)  &=&  a^{-1}\int_{T_n
}^{\infty}\int_{T_n }^{\infty}
\psi_\lambda\Big(\frac{t-b}{a}
\Big)\psi_\lambda\Big(\frac{t'-b}{a} \Big)\E \big (
X(t)X(t')\big )dtdt'
\\
 &\le&  C_f\, (a\, \lambda_n)^{-1}\,\Big (\int_{T_n }^{\infty}
\left|\psi\Big(\frac{t-b}{a \lambda_n} \Big)\right|\,\big(1+|t|\big)\,dt
\Big)^2\\
 &\le&  C_f\, (a^3 \lambda^3_n)\,\Big (\int_{ T_n^\rho/a \lambda_n }^{\infty}
|\psi(u )|\,u\,du
\Big)^2\\
 &\le& \frac 1 9 \, C_f\, C_\psi \, (a^3 \lambda^3_n)\,\Big (\Big [u^{-3} \Big ]_{T_n^\rho/a\lambda_n }^{\infty}\Big )^2\leq \frac 1 9 \, C_f\, C_\psi \,a^9\,  \lambda^9_n \, T_n^{-6\rho}.
\ean
Therefore the CLT holds when
$\ds \lambda_n^9 \, (n \, \delta_n)^{1-6\rho}  \limiten  0$, {\it i.e.} $\ds \frac {n \delta_n} {\lambda_n^2} \limiten \infty$ since $\rho>3/4$.
\end{dem}
\end{document}